\documentclass[a4paper,14]{elsarticle}
\usepackage[left=3cm,top=3cm,right=3cm,bottom=3cm]{geometry}
\usepackage[utf8]{inputenc}
\usepackage{times} %PostScript Times Roman instead of Computer Times Roman.
\usepackage{latexsym,amsfonts,amssymb,amsthm,amsmath,amscd,mathrsfs}
\usepackage{xargs}   % Use more than one optional parameter in a new commands
\usepackage{xspace}
\usepackage[shortlabels]{enumitem}
\usepackage{float}
\usepackage{multirow}
\usepackage[tight]{subfigure}
\usepackage{hyperref}
\usepackage{graphicx,tikz,times,pgffor}
\usepackage{pgfplots}
\usetikzlibrary{calc,through,backgrounds,patterns}
\usepackage{tabulary}
\usepackage{tabularx}
\usepackage{changepage}
\usepackage{soul}

\pgfplotsset{compat=newest} 
\pgfplotsset{plot coordinates/math parser=false} 
\newlength\figureheight 
\newlength\figurewidth 

\usepackage{xpatch}
\usepackage{algpseudocode}
\usepackage{algorithm}
\newcommand{\x} {\mathrm{x}}

\newcommand{\setX}{\mathcal{X}}

\newcommand{\setU}{\mathcal{U}}

\newcommand{\RN}{\mathcal{R}_0}
\newcommand{\Rn}{\mathcal{R}}
\newcommand{\Ko}{\mathcal{K}_0}

\newcommand{\Uc}{\mathcal{U}_c}

\newcommand{\I}{\mathbb{I}}
\newcommand{\X}{\mathbb{X}}

\newcommand{\R}{\mathbb{R}}

\newcommand {\bmat} {\left[\begin{array} }
\newcommand {\emat} {\end{array}\right]}

\newtheorem {theorem}{Theorem}
\newtheorem {defi}{Definition}
\newtheorem{rem}{Remark}
\newtheorem{lem}{Lemma}
\newtheorem {propt}{Property}

\newenvironment{pf}{\noindent \textbf{Proof:}}{}
\newenvironment{pf*}[1][Proof:]{\noindent \textbf{#1} }{}

\usepackage[shadow,colorinlistoftodos,prependcaption,spanish,textsize=footnotesize]{todonotes}
\newcommandx{\falta}[2][1=]{\todo[linecolor=red,backgroundcolor=red!25,bordercolor=red,#1]{#2}}
\newcommandx{\completar}[2][1=]{\todo[linecolor=blue,backgroundcolor=blue!25,bordercolor=blue,#1]{#2}}
\newcommandx{\chequear}[2][1=]{\todo[linecolor=OliveGreen,backgroundcolor=OliveGreen!25,bordercolor=OliveGreen,#1]{#2}}
\newcommandx{\improve}[2][1=]{\todo[linecolor=Plum,backgroundcolor=Plum!25,bordercolor=Plum,#1]{#2}}
\newcommandx{\explain}[2][1=]{\todo[linecolor=lime,backgroundcolor=lime!25,bordercolor=lime,#1]{#2}}
\newcommandx{\explicar}[2][1=]{\todo[linecolor=Maroon,backgroundcolor=Maroon!25,bordercolor=Maroon,#1]{#2}}
\definecolor{ao_english}{rgb}{0.0, 0.5, 0.0}
\definecolor{cagreen}{rgb}{0.12, 0.3, 0.17}
\definecolor{cgreen}{rgb}{0.0, 0.42, 0.24}
\definecolor{calpolypomonagreen}{rgb}{0.12, 0.3, 0.17}
\definecolor{armygreen}{rgb}{0.29, 0.33, 0.13}
\definecolor{myred}{rgb}{0.76, 0.23, 0.13}
\definecolor{bostonuniversityred}{rgb}{0.8, 0.0, 0.0}
\definecolor{blue(pigment)}{rgb}{0.2, 0.2, 0.6}
\definecolor{ceruleanblue}{rgb}{0.16, 0.32, 0.75}
\definecolor{darkblue}{rgb}{0.0, 0.0, 0.55}

\numberwithin{equation}{section}
\numberwithin{theorem}{section}

\allowdisplaybreaks
\begin{document}
\begin{frontmatter}

\title{Characterization of SARS-CoV-2 Dynamics in the Host}

\author[First]{Pablo Abuin} 
%\author[First]{Juan Sereno} 
\author[First]{Alejandro Anderson} 
\author[Second]{Antonio Ferramosca} 
%\author[Third]{Gustavo Hernandez-Mejia}
\author[Thirda,Third]{Esteban A. Hernandez-Vargas}
\author[First]{Alejandro H. Gonzalez} 

\address[First]{Institute of Technological Development for the Chemical Industry (INTEC), CONICET-UNL, Santa Fe, Argentina}
\address[Second]{CONICET - Universidad Teconologica Nacional, Facultad Regional Reconquista, (Santa Fe, Argentina).}
\address[Thirda]{Instituto de Matemáticas, Universidad Nacional Autonoma de Mexico, Boulevard Juriquilla 3001, Querétaro, Qro., 76230, México}
\address[Third]{Frankfurt Institute for Advanced Studies, 60438 Frankfurt am Main, Germany.\\
Corresponding author: alejgon@santafe-conicet.gov.ar, vargas@fias.uni-frankfurt.de
}

%%%%%%%%%%%%%%%%%%%%%%%%%%%%%%%%%%%%%%%%%%%%%%%%%%%%%%%%%%%%%%%%%%%%%%%%%%%%%%%%
\begin{abstract}
While many epidemiological models have being proposed to understand and handle COVID-19, too little has been invested to understand how the virus replicates in the human body and potential antiviral can be used to control the replication cycle. 
In this work, using a control theoretical approach, validated mathematical models of SARS-CoV-2 in humans are properly characterized. A complete analysis of the main dynamic characteristic is developed based on the
reproduction number. 
The equilibrium regions of the system are fully characterized, and the stability of such a regions, formally established. 
Mathematical analysis highlights critical conditions to decrease monotonically SARS-CoV-2 in the host, such conditions are relevant to tailor future antiviral treatments. 
Simulation results show the potential benefits of the aforementioned system characterization.
\end{abstract}
\end{frontmatter}

%%%%%%%%%%%%%%%%%%%%%%%%%%%%%%%%%%%%%%%%%%%%%%%%%%%
\section{Introduction}\label{sec:intro}
%%%%%%%%%%%%%%%%%%%%%%%%%%%%%%%%%%%%%%%%%%%%%%%%%%%

By December 2019, an outbreak of cases of pneumonia of unknown  etiology was reported in Wuhan, Hubei province, China \cite{luoutbreak}. On January 7, a novel betacoronavirus was identified as the etiological agent by the Chinese Center of Disease Control and Prevention (CDC), and subsequently named as Severe Acute Respiratory Syndrome Coronavirus 2 (SARS-CoV-2) \cite{gorbalenya2020severe}. On February 11, the World Health Organization (WHO) named the disease as Coronavirus disease 2019 (COVID-19) \cite{WHOCovid19}. Although, prevention and control measures were implemented rapidly, from the early stages in Whuan and other key areas of Hubei \cite{WHOChinaReport}, the first reporting of cases outside of China, 2 in Thailand and 1 in Japon \cite{WHOCovid19Timeline}, showed that the disease was starting to spread around the world. On March 11, with more that 111.800 cases in 114 countries, and 4921 fatality cases, the COVID-19 is declared pandemic by the WHO \cite{WHOCovid19Timeline}.  
So far, with more than 7.000.000 total cases confirmed in 213 countries and territories \cite{WHOCovid19SituationReport,Covid19Hopkins}, and a estimated case-fatality rate (CFR) of 5.7\% (H1N1 pandemic, CFR$<$1\%)\cite{WHOCovid19FatalityRate}, the potential risks associated with this disease are evident.

Facing this situation, and taking into account the nonexistence of vaccines or specific therapeutically treatments, preventive measures such as social and physical distancing, hand washing, cleaning and disinfection of surfaces and use of face masks, among others, have been implemented in order to decrease the transmission of the virus, which is spread mainly from person-to-person through respiratory droplets produced when an infected person coughs, sneezes or talks \cite{CDCCovid19Spread}. Furthermore, this infection prevention and control measures, helps avoid healthcare systems from becoming overwhelmed.

Many epidemiological mathematical models \cite{giordano2020sidarthe,acuna2020sars,read2020novel}  have been proposed to predict the spread of the disease and evaluate the potential impact of infection prevention and control measures in outbreak management \cite{anderson2020will}. 
However, mathematical models at within-host level that could be useful to understand SARS-CoV-2 replication cycle and interaction with immune system as well as pharmacological effect of potential drug therapies \cite{liu2020research, mitja2020use} are needed. 
So far, there are approximately 109 trials (including those not yet recruiting, recruiting, active, or completed) to asses pharmacological therapy for the treatment of COVID-19 in adult patients\cite{sanders2020pharmacologic}, including antiviral drugs (i.e. Hydroxychloroquine, Remdesivir, Favipiravir, Lopinavir/Ritonavir, Ribavirin), immunomodulatory agents (i.e. Tocilizumab) and immunoglobulin therapy, among others.

Recently, Hernandez-Vargas et. al. \cite{vargas2020host} proposed different within-host mathematical models (2 based on target cell-limited model, with and without latent phase, and 1 considering immune response) for 9 infected patients with COVID-19. Numerical results in \cite{vargas2020host} showed a mean infecting time between susceptible cells of 30 days (about 3 times slower than Ebola and 60 times slower than influenza), which could explain the slow recuperation rate (12-22 days post symptom onset, pso) showed in COVID-19 infected patients. Furthermore, they informed within-host reproductive number values consistent to influenza infection (1.7-5.35).

Although models in \cite{vargas2020host} have been fitted to COVID-19 patients data, a control theoretical approach is needed to characterize the model dynamics. Even when the equilibrium states are known, a formal stability analysis is needed to properly understand the model behavior and, mainly, to properly design efficient control strategies. Note that the target cell model has been employed previously taking into account pharmacodynamic (PD) and pharmacokinetics (PK) models of antiviral therapies \cite{hernandez2019passivity,boianelli2016oseltamivir}, and this can be potentially done also for COVID-19.

In this context, the main contribution of this article is twofold. First, a full characterization of equilibrium and stability proprieties is performed for the COVID-19 target cell-limited model \cite{vargas2020host}. Then, formal properties concerning the state variables behavior before convergence - including an analysis of the virus peak times - are given. 
A key aspect in the target cell model for acute infections shows some particularities such as 
it has a minimal non-punctual stable equilibrium set, whose stability does not depends on the reproduction number. On the other side, assuming a basic reproduction number greater than $1$, the virus would not be cleared before the target cells goes under a given critical value, which is independent of the initial conditions.

After the introduction given in Section \ref{sec:intro} the article is organized as follows. 
Section \ref{sec:infectmod} presents the general in-host target cell-limited  model used to represent SARS-CoV-2 infection dynamic. 
Section \ref{sec:assest} characterizes the equilibrium sets of the system, and establishes their formal asymptotic stability, by
stating both, the attractivity of the equilibrium set in a given domain, and its $\epsilon-\delta$ (Lyapunov) local stability.
Then, in Section \ref{sec:dyn_char}, some dynamical properties of the system are stated, concerning the values of the
states at the infection time $t=0$. 
In Section \ref{sec:pac} the general model for the SARS-CoV-2 infection is identified according to patient data, and the general characteristics of the infection are analyzed. Finally, Section \ref{sec:conc} gives the conclusion of the work, while several mathematical formalism - necessary to support the results of Sections \ref{sec:assest} and \ref{sec:dyn_char} - are given in Appendices \ref{sec:app1}, \ref{sec:app2} and \ref{sec:app3}.

%%%%%%%%%%%%%%%%%%%%%%%%%%%%%%%%%%%%%%%%%%%%%%%%%%%%%%%%%%%%%%%%%%%%%%%%%%%%%%%%%%%%%%%%%%%%%%%%
\subsection{Notation}\label{sec:not}
%%%%%%%%%%%%%%%%%%%%%%%%%%%%%%%%%%%%%%%%%%%%%%%%%%%%%%%%%%%%%%%%%%%%%%%%%%%%%%%%%%%%%%%%%%%%%%%%

$\R$ and $\I$ denote the real and integer numbers, respectively. The real vector space of dimension $n$ is denoted as $\R^n$. 
$\R^n_{\geq 0}$ represents the vectors of dimension $n$ whose components are equal or greater than zero. The distance from a point $x \in \R^n$ to a set $\setX \subset \R^n$ is defined by $\|x\|_{\setX}:= \inf_{z \in \setX} \|x-z\|_2$, where $\|\cdot\|_2$ denotes the norm-$2$. The open ball of radius $\epsilon$ around a point $x \in \R^n$, with respect to set $\setX$, is defined as $\mathbb{B}_{\epsilon}(x):= \{z\in \setX: \|x-z\|_2 < \epsilon \}$. Let us consider the real function $y(z)=ze^z$, then, the so-called Lambert function is defined as the inverse of $y(\cdot)$, i.e., $W(z):=f^{-1}(z)$ in such a way that $W(f(z))=z$.

%%%%%%%%%%%%%%%%%%%%%%%%%%%%%%%%%%%%%%%%%%%%%%%%%%%%%%%%%%%%%%%%%%%%
\section{SARS-CoV-2 Within-Host Mathematical Model}\label{sec:infectmod}
%%%%%%%%%%%%%%%%%%%%%%%%%%%%%%%%%%%%%%%%%%%%%%%%%%%%%%%%%%%%%%%%%%

Although incomplete by definition, mathematical models of in-host virus dynamic improves the understanding of the interactions that
govern infections and, more important, permits the human intervention to moderate their effects \cite{hernandez2019modeling}. Basic in-host infection dynamic models usually include the susceptible cells, infected cells, and the pathogen particles \cite{ciupe2017host}. 
Among the most used mathematical models, the target cell-limited model has been employed to represent and control HIV infection \cite{perelson1993dynamics,legrand2003vivo,perelson2013modeling}, 
influenza \cite{larson1976influenza,baccam2006kinetics,smith2011influenza,hernandez2019passivity}, Ebola \cite{nguyen2015ebola}, dengue \cite{nikin2015role,nikin2018modelling} among others.
%\textcolor{ao_english} {A main distinction of target cell models can be done according to infection period length \cite{ciupe2017host}. This way, we have chronic (i.e. long-lived and persistent infections in comparison with body cells regeneration rate) and acute (short-lived) infections, which are modeled quite differently.  Chronic infection models (first used to model HIV in-host viral kinetics) are indeed a generalization of acute infection ones, that adds the body production and death rates of healthy cells. This is so, because the infection length is comparable with these rates \cite{ciupe2017host}}. \todo{Ale G: Pablo, please improve this first paragraph.}

In this work, we consider the mathematical model proposed by Hernandez-Vargas \cite{vargas2020host} given by the following set of differential equations (ODEs) :
\begin{subequations}\label{eq:SysOrigAcut}
	\begin{align}
	&\dot{U}(t)   =  -\beta U(t) V(t),~~~~~~~ U(0) = U_0, \label{eq:SysOrigAcut_healty}\\
	&\dot{I}(t) =  \beta U(t) V(t) - \delta I(t),~~~~~~~I(0) = I_0=0,\label{eq:SysOrigAcut_inf}\\
	&\dot{V}(t)   = pI(t) - c V(t),~~~~~~~ V(0) = V_0, \label{eq:SysOrigAcut_vir}
	\end{align}
\end{subequations}
where $U$ $[cell]$, $I$ $[cell]$ and $V$ $[copies/mL]$ represent the susceptible cells, the infected cells, and the virus load, respectively.
The parameter $\beta$ $[(copies/mL)^{-1}day^{-1}]$ is the infection rate of susceptible cells by the virus. 
$\delta$ $[day^{-1}]$ is the death rate of $I$.
$p$ $[(copies/mL)day^{-1}cell^{-1}]$ is the replication rate of free virus from infected cell $I$.
$c$ $[day^{-1}]$ is the degradation (or clearance) rate of virus $V$. 
The effects of immune responses are not explicitly described in this model, but they are implicitly included in the death rate of infected cells ($\delta$) and the clearance rate of virus ($c$) \cite{baccam2006kinetics}.

The model \eqref{eq:SysOrigAcut} is positive, which means that $U(t) \geq 0$, $I(t) \geq 0$ and $V(t) \geq 0$, for all $t\geq0$.
If we denote $x(t):=(U(t),I(t),V(t))$, then the states are constrained to belong to:
\begin{eqnarray}\label{eq:setX}
\mathbb X :=\{x\in \mathbb R^3_{\geq0}\}.
\end{eqnarray}

Another meaningful set is the one consisting in all the states in $\mathbb{X}$ with strictly positive amount of virus and 
susceptible cells, i.e.,

\begin{eqnarray}\label{eq:setXcal}
\setX :=\{x\in \mathbb X : U>0,~V > 0\}.
\end{eqnarray}
Note that the set $\setX$ is an open set.

The initial conditions of (\ref{eq:SysOrigAcut}) must be carefully established in order %for the mathematical model 
to properly represent the host body evolution from the beginning of the infection. So, it is assumed that the system is at a healthy steady state before the infection time $t=0$, i.e., $V(t) = 0$, $I(t) = 0$, and $U(t)= U_0$, for $t<0$. 
At time $t=0$, a small quantity of virions enters to the host body and, so, a discontinuity occurs in $V(t)$. Indeed, $V(t)$ jumps from $0$ to a small positive value $V_0$ at $t_0=0$ (formally, $V(t)$ has a discontinuity of the first kind at $t_0$, i.e., $\lim_{t\to0^-} V(t)=0$ while $\lim_{t\to0^+} V(t)=V_0>0$. 
The same scenario arises, for instance, when an antiviral treatment affects either  parameter (say $p$ or $\beta$).% starts at some time $t_{tr}>0$, and the \textcolor{ao_english}{pharmacokinetic} is neglected (the drug evolves in the body sufficiently fast compared \textst{to the virus} \textcolor{ao_english}{with the antiviral effect}). %on the viral load}). 
The jump of $p$ or $\beta$ can be considered as a discontinuity of the first kind.
In any case, for the time after the discontinuity, the virus may spread or be cleared in the body, depending on its infection effectiveness. To properly determine what such a spread means, the following (mathematical) definition is given
\begin{defi}[Spreadability of the virus in the host body]\label{defi:spread}
	Consider system \eqref{eq:SysOrigAcut}, constrained by the positive set $\X$, at some time $t_0$, with
	$U(t_0) >0$, $I(t_0) \geq 0$ and $V(t_0) > 0$ (i.e., $x(t_0)=(U(t_0),I(t_0),V(t_0)) \in \setX$). Then, it is said that 
	the virus spreads (in some degree) in the body host for $t>t_0$ if there exists at least one $t^*>t_0$
	such that $\dot{V}(t^*)>0$. 
\end{defi}
The latter definition states that the virus spreads in the body host if $V(t)$ has at least one local maximum. On the other hand, the virus does not spread if $V(t)$ is strictly decreasing for all $t>t_0$, which means that $V(t)$ has neither local minima nor local maxima. As it will be stated later on (Property \ref{prop:inftycond}), $\lim_{t\to\infty} V(t)=0$ for system \eqref{eq:SysOrigAcut}, independently of the fact that the virus reaches or not a maximum (this is a key difference between acute and chronic infection models \cite{hernandez2019modeling,ciupe2017host}). % Indeed, what can produce an infection and a disease is such a pick of virus, if severe. So, this is the reason why the spreadability defined in Definition \ref{defi:spread} concerns to virus maxima and not to virus level. 

The infection severity could be related with the virus spreadability established in Definition \ref{defi:spread}. Liu et.al. \cite{liu2020viral} have shown that patients with severe COVID-19 tend to have a high viral load and a long virus shedding period. The mean viral load of severe cases was around 60 times higher than that of mild cases, suggesting that higher viral loads might be associated with severe clinical outcomes. Furthermore, they found that the viral load of severe cases remained significantly higher for the first 12 days after the appearance of the symptoms than those of corresponding mild cases. Mild cases were also found to have an early viral clearance, with 90\% of these patients repeatedly testing negative on reverse transcription polymerase chain reaction (RT-PCR) by day 10 post symptoms onset (pso). By contrast, all severe cases still tested positive at or beyond day 10 pso. In addition, Zheng et.al. \cite{zheng2020viral} from an study with 96 SARS-CoV-2 positive patients (22 with mild disease and 74 with severe disease) reported a longer duration of SARS-CoV-2 in lower respiratory samples of severe patients, such as, for patients with severe disease (21 days, 14-30 days) was significantly longer than in patients with mild disease (14 days, 10-21 days; p=0.04). Moreover, higher viral loads were detected in respiratory samples, although no differences were found in stool and serum samples.  Although, these findings suggest that reducing the viral load through clinical means and strengthening management should help to prevent the spread of the virus, they are preliminary and it remains controversial whether virus persistence is necessary to drive the dysfunctional immune response characteristic of COVID-19 patients \cite{tay2020trinity}. 

\begin{rem}
	Note that the virus spreadability may or may not cause a severe infection (a disease that eventually cause the host death) depending on how much time the virus is above a given value.
\end{rem}

To properly establish conditions under which the virus does not spread for $t>0$ (i.e., after the infection time $t=0$)
the so-called basic reproduction number within-host is defined next.
\begin{defi}\label{defi:repnumb}
	The basic reproduction number within-host $\Rn$ is defined as the number of infected cells (or virus particles) that are produced by
	one infected cell (or virus particle), at a given time. Its mathematical expression is given by:
	\begin{eqnarray}\label{eq:Rn}
		\Rn(t) := U(t)\frac{\beta p}{c \delta} .
	\end{eqnarray}
	Particularly, for $t=0$, this number describes the number of infected cells produced by  one infected cell, when a small amount of virus, $V_0$, is introduced into a healthy stationary population of uninfected target cells, $U_0$,
	\begin{eqnarray}\label{eq:RN}
		\RN := U_0\frac{\beta p}{c \delta} .
	\end{eqnarray}
\end{defi}
A discussion about the way this value is obtained is given in Appendix \ref{sec:app2}. The relation between the basic reproduction
number at the infection time ($\RN$) and the virus spreadiblity is stated in the next theorem. 
\begin{theorem}\label{the:repnumber}
	Consider system \eqref{eq:SysOrigAcut}, constrained by the positive set $\X$, at the beginning of the infection, i.e.,
	$U(0) = U_0 >0$, $I(0)= 0$ and $V(0)=V_0>0$ (i.e., $x(0)=(U(0),I(0),V(0)) \in \setX$). Then, a sufficient condition (not necessary) 
	for the virus not to spread in the host body is given by $\RN<1$.
\end{theorem}
\begin{pf}
	It is easy to see that the initial conditions are such that $\Rn_V(0):=\frac{pI(0)}{cV(0)}=0<1$. So, considering that
	$\RN <1$ by hypothesis, the result corresponds to Theorem \ref{theo:key}. (i)., when $t_0=0$.  $\square$
\end{pf}
%\vspace{0.5cm}
\begin{rem}
	According to Theorem \ref{theo:key}, in Appendix \ref{sec:app3}, there exists a positive value $\alpha(0)>0$ (which is an implicit function of 
	the initial conditions and the parameters, and it can be computed numerically) such that condition $\RN<1+\alpha(0)$ is 
	a necessary and sufficient condition for the virus not to spread in the host body, after the infection time $t=0$.
	This means that for $\RN>1+\alpha(0)$ the virus spreads in the host, as it is usually the case in real infections. 
	In any case, the value of $\alpha(0)$ is generally close to zero (when $\beta$ is small), so to consider $1$ as a threshold 
	for the virus to spread seems to be a reasonable approximation.
\end{rem}

Before proceeding with a full dynamic analysis of system \eqref{eq:SysOrigAcut}, let us define first the so-called
critical value of the susceptible cells, which shows to be an useful threshold value to properly understand the spread of the virus.   
\begin{defi}
	The critical value for $U$, $\Uc$, is defined as
\begin{eqnarray}\label{eq:ucrit}
	\Uc := \frac{c \delta}{p \beta},
\end{eqnarray}
which, for fixed system parameters $\beta$, $p$, $\delta$ and $c$, is a constant.
\end{defi}
Note that $U(t) < \Uc$ if and only if $\Rn(t) <1$, for every $t\geq0$.

%%%%%%%%%%%%%%%%%%%%%%%%%%%%%%%%%%%%%%%%%%%%%%%%%%%%%
\subsection{Equilibrium set characterization}\label{sec:equilset}
%%%%%%%%%%%%%%%%%%%%%%%%%%%%%%%%%%%%%%%%%%%%%%%%%%%%%

By equaling $\dot U$, $\dot I$ and $\dot V$ to zero, in \eqref{eq:SysOrigAcut}, it can be shown 
that the system only has healthy equilibria of the form $x_s=(U_s,0,0)$, with $U_s$ being an arbitrary 
positive value, i.e., $U_s \in [0,\infty)$. Thus, there is only one equilibrium set, which is the healthy one, and it is defined by 
\begin{eqnarray}\label{eq:equilset}
\setX_s := \{(U,I,V)\in \mathbb{R}^3 : U \in [0,\infty),~I=0,~V=0\}.
\end{eqnarray}
To have a first look on the stability of the equilibrium points in $\setX_s$, system \eqref{eq:SysOrigAcut} can be linearized at a 
general state $x_s \in \setX_s$.
By simplifying \eqref{eq:SysOrigAcut} we have:
\begin{eqnarray*}
	\dot{U}   &=&  f(U,I,V), \\
	\dot{I}   &=&  g(U,I,V),\\
	\dot{V}   &=&  h(U,I,V).
\end{eqnarray*}
Then, the Jacobian matrix is given by 

\begin{eqnarray*}
	J = \left( 
	\begin{array}{ccc}
		\displaystyle \frac{\partial f}{\partial U}   & \displaystyle\frac{\partial f}{\partial I} &   \displaystyle\frac{\partial f}{\partial V}   \\
		\displaystyle\frac{\partial g}{\partial U}  & \displaystyle\frac{\partial g}{\partial I} &   \displaystyle\frac{\partial g}{\partial V}\\
		\displaystyle\frac{\partial h}{\partial U}   & \displaystyle\frac{\partial h}{\partial I} &   \displaystyle\frac{\partial h}{\partial V}
	\end{array}\right) =
	\left( 
	\begin{array}{ccc}
		-\beta V  & 0 &   -\beta U   \\
		\beta V   &  -\delta & \beta U  \\
		0    & p & -c   
	\end{array}\right).
\end{eqnarray*}
And the Jacobian evaluated at any point $x_s \in \setX_s$ reads
\begin{eqnarray*}
	A_s = \left( 
	\begin{array}{ccc}
		0   &    0 &  -\beta U_s   \\
		0   &  -\delta & \beta U_s  \\
		0   & p & -c   
	\end{array}\right).
\end{eqnarray*}
with $U_s \in [0,\infty)$. Then, the eigenvalues $(\lambda_1,\lambda_2,\lambda_3)$ are computed as the solution to $Det(A_s-\lambda I)=0$,
being the matrix $A_s-\lambda I$ given by
\begin{eqnarray*}
	A_s-\lambda I = \left( 
	\begin{array}{ccc}
		-\lambda   &    0 &  -\beta U_s   \\
		0   &  (-\delta-\lambda) & \beta U_s  \\
		0   & p & (-c-\lambda)  
	\end{array}\right).
\end{eqnarray*}
Then, considering that $Det(A-\lambda I) = \lambda[-\lambda^2 - (c+\delta) \lambda + (\beta U_sp-c\delta)]$, 
condition $Det(A-\lambda I)=0$ is given by
\begin{eqnarray*}
	\lambda[-\lambda^2 - (c+\delta) \lambda + (\beta U_sp-c\delta)]=0.
\end{eqnarray*}
The first eigenvalue is trivially given by $\lambda_1=0$. The other two, are given by:
\begin{eqnarray*}
	\lambda_{2,3}= \frac{(c+\delta) \pm \sqrt{(c+\delta)^2 + 4(\beta U_s p - c \delta)}}{-2}.
\end{eqnarray*}
To analyze the eigenvalues qualitatively, note that for $U_s = \Uc$ it is
\begin{eqnarray*}
	\lambda_{2,3} &=& \frac{(c+\delta) \pm \sqrt{(c+\delta)^2 + 4(\beta p \frac{c\delta}{p \beta} - c \delta)}}{-2}\\
	&=& \frac{(c+\delta) \pm \sqrt{(c+\delta)^2 + 4(c\delta - c \delta)}}{-2}\\
	&=& \frac{(c+\delta) \pm (c+\delta)}{-2},
\end{eqnarray*}
which means that $\lambda_2=0$ and $\lambda_3 = -(c+\delta) < 0$ (given that $c,\delta>0$). Furthermore, 
$\lambda_2 <0$ and $\lambda_3<0$ for $U_s < \Uc$; and $\lambda_2 >0$ and $\lambda_3<0$ for $U_s > \Uc$.
Given that the maximum eigenvalue is the one dominating the stability behavior of the equilibrium under consideration,
it is possible to infer how the system behaves
near some segments of $\setX_s$. The first intuition is that the equilibrium set 
\begin{eqnarray}\label{ec:setXs1}
\setX_s^1 := \{(U,I,V)\in \mathbb{R}^3 : U \in [0,\Uc),~I=0,~V=0\}
\end{eqnarray}
is stable, and that the equilibrium set 
\begin{eqnarray}\label{ec:setXs2}
\setX_s^2 := \{(U,I,V)\in \mathbb{R}^3 : U \in [\Uc,\infty),~I=0,~V=0\}
\end{eqnarray}
is unstable. These are just intuitions, given that one of the eigenvalues of the linearization system is null and so the linear approximation cannot be 
used to fully determine the stability of the nonlinear system (Theorem of Hartman-Grobman \cite{hartman1982ordinary,perko2013differential}). 
To formally prove the asymptotic stability of $\setX_s^1$ in a given domain, it is necessary to prove its global attractivity (in such domain) and local 
$\epsilon$-$\delta$ stability.

%%%%%%%%%%%%%%%%%%%%%%%%%%%%%%%%%%%%%%%%%%%%%%%%%%%%%
\section{Asymptotic stability of the equilibrium sets}\label{sec:assest}
%%%%%%%%%%%%%%%%%%%%%%%%%%%%%%%%%%%%%%%%%%%%%%%%%%%%%

A key point to analyze the general asymptotic stability (AS) of system \eqref{eq:SysOrigAcut} is to consider stability of the complete 
equilibrium sets $\setX_s^1$ and $\setX_s^2$, and not of the single points inside them (as defined in Definitions \ref{def:attrac_set}, \ref{def:eps_del_stab}
and \ref{def:AS}, in Appendix \ref{sec:app1}). As it is shown in the next subsections, there is no 
single AS equilibrium points in this system, although there is an AS equilibrium set (i.e., $\setX_s^1$). 

As stated in Definition \ref{def:AS}, in Appendix \ref{sec:app1}, the AS of $\setX_s^1$ requires both, attractivity 
and $\epsilon-\delta$ stability, which are stated in the next two subsections, respectively. Then, in Subsection \ref{sec:ASproof}
the AS theorem is formally stated.

%%%%%%%%%%%%%%%%%%%%%%%%%%%%%%%%%%%%%%%%%%%%%%%%%%%%%
\subsection{Attractivity of set $\setX_s^1$ in $\setX$}\label{sec:attrac}
%%%%%%%%%%%%%%%%%%%%%%%%%%%%%%%%%%%%%%%%%%%%%%%%%%%%%

Before proceeding with the formal theorems of the atractivity of $\setX_s^1$, let us consider the 
following key property of system \eqref{eq:SysOrigAcut} concerning the atractivity of $\setX_s$.
\begin{propt}[Atractivity of $\setX_s$]\label{prop:inftycond}
	Consider system \eqref{eq:SysOrigAcut} constrained by the positive set $\mathbb X$, at some arbitrary time $t_0$, with
	$U(t_0) >0$, $I(t_0) \geq 0$ and $V(t_0) > 0$ (i.e., $x(t_0)=(U(t_0),I(t_0),V(t_0)) \in \setX$).
	Then, $U_{\infty}:=\lim_{t \rightarrow \infty}U(t)$ is a constant value smaller than $U(t_0)$, $I_{\infty}:=\lim_{t \rightarrow \infty}I(t) = 0$ and $V_{\infty}:=\lim_{t \rightarrow \infty}V(t) = 0$, which means that $x(t)=(U(t),I(t),V(t))$ tends to some state in $\setX_s$.
\end{propt}
\begin{pf}
	Since $\dot{U}(t) \leq 0$ for all $t \geq 0$ and all $(U(t_0),I(t_0),V(t_0)) \in \setX$, by (\ref{eq:SysOrigAcut_healty}) $U(t)$ is a decreasing function (no oscillation can occur).
	Since $U(t_0) > 0$ and $V(t_0)>0$, then $U_{\infty}=\lim_{t \rightarrow \infty}U(t)$ is a constant value in $[0,U(t_0))$. Given that $U(t)$ converges to a finite fixed value, 
	then $\dot{U}(t)=0$ as $t \rightarrow \infty$, by (\ref{eq:SysOrigAcut_healty}). 
	This implies, by the same equation (\ref{eq:SysOrigAcut_healty}), that $U(t)V(t) =0$ as $t \rightarrow \infty$, and so, from equation (\ref{eq:SysOrigAcut_inf}),
	that $\dot I(t) = -\delta I(t)$ as $t \to \infty$, whose solution asymptotically goes to zero. Then, $I_{\infty} = \lim_{t \rightarrow \infty} I(t) = 0$.
	Finally, by equation (\ref{eq:SysOrigAcut_vir})), $\dot V(t) = -\delta V(t)$ as $t \rightarrow \infty$, whose solution asymptotically goes to zero.
	Then $V_{\infty}=\lim_{t \rightarrow \infty} V(t) =0$, which completes the proof. $\square$
\end{pf}
\vspace{0.2cm}

Property \ref{prop:inftycond} states that $\setX_s$ is an attractive set for system \eqref{eq:SysOrigAcut}, in $\setX$, but not the 
smallest attractive set. Now, conditions are given to show that the smallest attractive set is given by $\setX_s^1$.
\begin{theorem}[Atractivity of $\setX_s^1$]\label{theo:attract}
	Consider system \eqref{eq:SysOrigAcut} constrained by the positive set $\mathbb X$.
	Then, the set $\setX_s^1$ defined in \eqref{ec:setXs1} is the smallest attractive set in $\setX$.
	Furthermore, $\setX_s^2$, defined in \eqref{ec:setXs2}, is not attractive.
\end{theorem}

\begin{pf}
	The proof is divided into two parts. First it is proved that $\setX_s^1$ is an attractive set, and then, that it is the smallest one.
	
	\textit{Attractivity of $\setX_s^1$:}
	The attractivity of $\setX_s$ in $\setX$ is already proved in Property \ref{prop:inftycond}. So, to prove the attractivity 
	of $\setX_s^1$ in $\setX$ (and to show that $\setX_s^2$ is not attractive) it remains to demostrate that $U_{\infty} \in [0,\Uc)$.
	From system \eqref{eq:SysOrigAcut}, by replacing (\ref{eq:SysOrigAcut_healty}) in (\ref{eq:SysOrigAcut_inf}), it follows that 
	$\dot{I}(t) =  \beta U(t) V(t) - \delta I(t) = -\dot{U}(t) - \delta I(t)$,
	%
%	\begin{eqnarray} \label{ec:dem1}
%	\dot{I}(t) =  \beta U(t) V(t) - \delta I(t) = -\dot{U}(t) - \delta I(t),
%	\end{eqnarray}
	%
	which implies that 
	\begin{eqnarray} \label{ec:dem2}
	I(t) = (-\frac{1}{\delta}) (\dot{I}(t) + \dot{U}(t)).
	\end{eqnarray}
	From (\ref{eq:SysOrigAcut_vir}) it follows that 
	\begin{eqnarray} \label{ec:dem3}
	V(t) =  \frac{1}{c}(p I(t) - \dot{V}(t)).
	\end{eqnarray}
	Then, replacing \eqref{ec:dem2} in \eqref{ec:dem3}, we have
	\begin{eqnarray} \label{ec:dem4}
	V(t) = [p (-\frac{1}{\delta}) (\dot{I}(t) + \dot{U}(t)) - \dot{V}(t)] \frac{1}{c}.
	\end{eqnarray}
	Finally, by substituting \eqref{ec:dem4} in (\ref{eq:SysOrigAcut_healty}), and multiplying by $1/U(t)$ both sides of the 
	equation (without loss of generality we assume that $U(t)\not=0$), it follows that 
	\begin{eqnarray} \label{ec:dem5}
	\frac{1}{U(t)} \dot{U}(t) = \frac{\beta p}{c \delta} \dot{U}(t) + \frac{\beta p}{c \delta} \dot{I}(t) + \frac{\beta}{c} \dot{V}(t) .
	\end{eqnarray}
	%
	%where $\frac{\beta p}{c \delta} = \Uc^{-1}$. 
	This latter equation can be integrated, for general initial conditions $U_0$, $I_0$ and $V_0$, as follows:
	\begin{eqnarray} \label{ec:dem6}
	\ln(\frac{U(t)}{U_0}) = \frac{\beta p}{c \delta} (U(t)-U_0) + \frac{\beta p}{c \delta} (I(t)-I_0) + \frac{\beta}{c} (V(t)-V_0).
	\end{eqnarray}
	Now, by defining $U_{\infty}:=\lim_{t \rightarrow \infty }U(t)$, $I_{\infty}:=\lim_{t \rightarrow \infty }I(t)$, $V_{\infty }:=\lim_{t \rightarrow \infty }V(t)$, and recalling from Property \ref{prop:inftycond} %(in Appendix \ref{sec:app3}) 
	that $I_{\infty }=V_{\infty }=0$, the latter equation for $t \rightarrow \infty$, reads
	\begin{eqnarray} \label{ec:dem7}
	\ln(\frac{U_{\infty}}{U_0}) &=& \frac{\beta p}{c \delta} (U_{\infty}-U_0) + \frac{\beta p}{c \delta} (I_{\infty}-I_0) + \frac{\beta}{c} (V_{\infty}-V_0) \nonumber \\
	&=& \frac{\beta p}{c \delta} U_{\infty} - \RN - \frac{\beta p}{c \delta} I_0 - \frac{\beta}{c} V_0 \nonumber \\
	&=& \frac{\beta p}{c \delta} U_{\infty} - \RN  + \Ko,
	\end{eqnarray}
	where $\RN:=\frac{\beta p}{c \delta}U_0$ (as it was defined in \eqref{eq:RN}) and
	\begin{eqnarray} \label{ec:dem8}
	\Ko := - \frac{\beta}{c} (\frac{p}{\delta} I_0 + V_0).
	\end{eqnarray}
	Note that $\RN$ is a function of $U_0$ while $\Ko$ is a function of $I_0$ and $V_0$, and, furthermore, $\RN>0$ and $\Ko<0$ for every 
	$x_0=(U_0,I_0,V_0) \in \setX$. Then, after some manipulation, \eqref{ec:dem7} reads
	\begin{eqnarray} \label{ec:dem8_5}
	-\frac{\beta p}{c \delta} U_{\infty} e^{-\frac{\beta p}{c \delta} U_{\infty}} = -\frac{\beta p}{c \delta}U_0 e^{-\RN} e^{\Ko} = -\RN e^{-\RN} e^{\Ko}.
	\end{eqnarray}
	Now, by denoting $z=z(\RN,\Ko):=-\RN e^{-\RN} e^{\Ko}$ and $y:=-\frac{\beta p}{c \delta} U_{\infty}$, the latter equation can be written as
	\begin{eqnarray} \label{ec:dem9}
	W(z)=y,
	\end{eqnarray}
	or, the same,
	\begin{eqnarray} \label{ec:dem10}
	W(-\RN e^{-\RN} e^{\Ko}) = -\frac{\beta p}{c \delta} U_{\infty},
	\end{eqnarray}
	where $W(\cdot)$ is a Lambert function. Figure \ref{fig:Lamb} shows the graph of such a function, where it can be seen that 
	it has two branches, denoted as $W_p$ and $W_m$. However, $W(\cdot)=W_p(\cdot)$ in this case, since $W_m \rightarrow -\infty$ 
 	for $z \rightarrow 0^-$, which has not biological sense (note that $U_{\infty}$ is a finite value in $[0,U_0)$). Besides, function $-1/e < z(\RN,\Ko) \leq 0$ for $\RN>0$ and $\Ko<0$ (Figure \ref{fig:FR1} shows a plot of function $z(\RN,\Ko)$ for negative values of $\Ko$ and positive values of $\RN$), and function $W_p$ maps $(-1/e,0]$ into $(-1,0]$, which implies that
	\begin{eqnarray} \label{ec:dem11}
	1 > -W(z(\RN,\Ko)) \geq 0,
	\end{eqnarray}
	for $\RN>0$ and $\Ko<0$. This way, by \eqref{ec:dem10}, it follows that 
	\begin{eqnarray} \label{ec:dem12}
	U_{\infty} &=&-\frac{c \delta}{\beta p}  W(-\RN e^{-\RN} e^{\Ko}) \nonumber\\
	           &=& - \Uc W(-\RN e^{-\RN} e^{\Ko}) \nonumber\\
	           &\in& [0,\Uc),
	\end{eqnarray}
	which completes the proof.
	\begin{figure}
		\centering
		\includegraphics[width=0.65\textwidth]{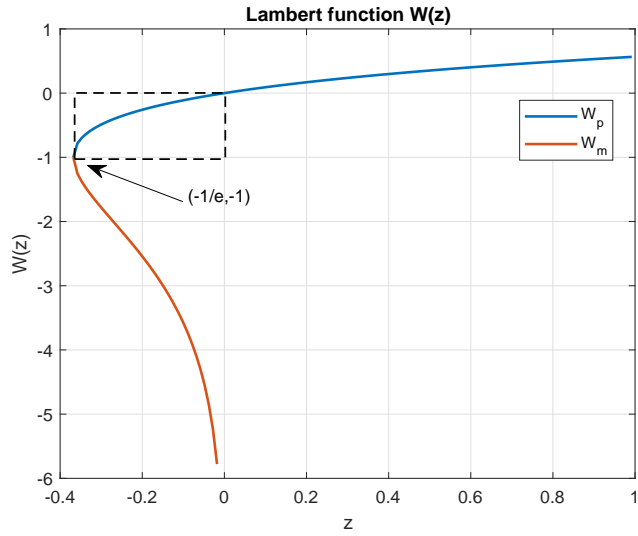}
		\caption{\small{Lambert function. $W(z)$ has two branches, denoted as $W_p$ (in blue) and $W_m$ (in red). Both branches are defined for
				$z \in [-1/e,0]$: however $\lim_{z \rightarrow 0^-} W_p=0$ while $\lim_{x \rightarrow 0^-} W_m = -\infty$, which means that 
				only the branch $W_p$ will be used in our analysis, as it is shown in the proof.}}
		\label{fig:Lamb}
	\end{figure}
	\begin{figure}
		\centering
		\includegraphics[width=0.65\textwidth]{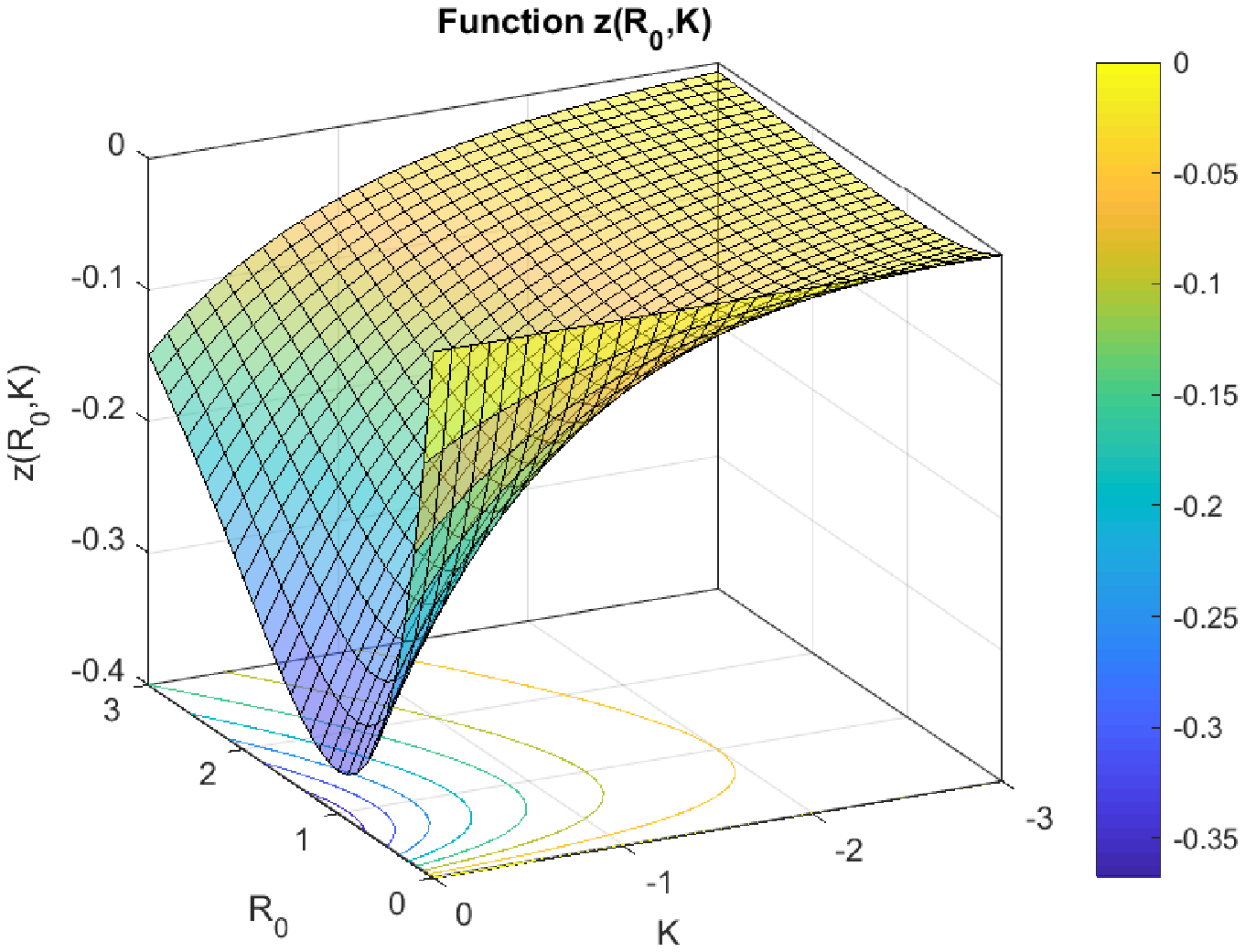}
		\caption{\small{Function $z(\RN,\Ko)$, for $\RN\geq 0$ and $\Ko \leq 0$.}}
		\label{fig:FR1}
	\end{figure}

	\textit{$\setX_s^1$ is the smallest attractive set:}
	It is clear, from the previous analysis, that any initial state $x_0 =(U_0,I_0,V_0)$ in $\setX$ converges to a state $x_{\infty}=(U_{\infty},0,0)$ with $U_{\infty} \in [0,\Uc)$. This means that every state $x_s \in \setX_s^2$ is not attractive in $\setX$ and so, neither the whole set $\setX_s^2$. Let us consider now a state $x_s \in \setX_s^1$ and an arbitrary small ball of radius $\epsilon > 0$, w.r.t. $\setX$, around it, $\mathbb{B}_{\epsilon}(x_s) \in \setX$. Take two arbitrary initial states $x_{0,1}=(U_{0,1},I_{0,1},V_{0,1})$ and $x_{0,2}=(U_{0,2},I_{0,2},V_{0,2})$ in $\mathbb{B}_{\epsilon}(x_s)$, such that 
 	$U_{0,1} \neq U_{0,2}$ and $V_{0,1} \neq V_{0,2}$. These two states converge, according to equation \eqref{ec:dem11},
	to $x_{\infty,1}=(U_{\infty,1},0,0)$ and $x_{\infty,2}=(U_{\infty,2},0,0)$, respectively. Given that function $z(R,K)$ is monotone \textcolor{black}{(injective)} in $\RN$ (and so in $U_0$) and $W(z)$ is monotone  \textcolor{black}{(injective)} in $z$, then $U_{\infty,1} \neq U_{\infty,2}$. This means that, although both initial states converge to some state in $\setX_s^1$, they necessarily converge to different points. Therefore neither single states $x_s \in \setX_s^1$ nor subsets of $\setX_s^1$ are attractive in $\setX$. So, $\setX_s^1$ is the smallest attractive set and the proof is concluded.
	$\square$
	
\end{pf}

\begin{rem}
	Note that $\setX_s^1$ and $\setX_s^2$ are in the closure of the open set $\setX$, which is not in $\setX$. In other words, what 
	Theorem \ref{theo:attract} shows is that any initial state in $\setX$ converges to a point onto the boundary of $\setX$ that does not belong to $\setX$. Furthermore note that, an initial state of the form $(U_0,0,0)$, $U_0>\Uc$, (i.e., a state in $\setX_s^2$) cannot be attracted by any set since it is - by definition - an equilibrium state (every state in $\setX_s^2$ will remains unmodified). This is the reason why it is not possible to consider the attractivity of $\setX_s^2$ in $\setX$.
\end{rem}

%%%%%%%%%%%%%%%%%%%%%%%%%%%%%%%%%%%%%%%%%%%%%%%%%%%%%%%%%%%%%%%%%%
\subsection{Local $\epsilon-\delta$ stability of $\setX_s^1$} \label{sec:epdelstabil}
%%%%%%%%%%%%%%%%%%%%%%%%%%%%%%%%%%%%%%%%%%%%%%%%%%%%%%%%%%%%%%%%%%

The next theorem states the formal Lyapunov (or $\epsilon-\delta$) stability of the equilibrium set $\setX_s^1$.

\begin{theorem}\label{theo:stab}
	Consider system \eqref{eq:SysOrigAcut} constrained by the positive set $\mathbb X$.
	Then, the equilibrium set $\setX_s^1$ defined in \eqref{ec:setXs1} is locally $\epsilon-\delta$ stable.
\end{theorem}
\begin{pf}
	Let us consider a particular equilbrium point $x_s := (U_s,0,0)$, with $U_s \in [0,\Uc)$ (i.e., $x_s \in \setX_s^1$). 
	Then a Lyapunov function candidate is given by (similar to one used in \cite{nangue2019global} for chronic infections)
	\begin{eqnarray} \label{ec:lya1}
	J(x) := U-U_s - U_s \ln(\frac{U}{U_s}) + I + \frac{\delta}{p} V.
	\end{eqnarray}
	This function is continuous in $\mathbb{X}$, is positive for all nonegative $x \neq x_s$ and $J(x_s)=0$. 
	Function $J$ evaluated at the solutions of system \eqref{eq:SysOrigAcut} reads:
	\begin{eqnarray} \label{ec:lya2}
	\dot{J}(x(t)) &=& \frac{\partial J}{\partial x} \dot{x}(t) =  \left[\frac{d J}{d U}~~\frac{d J}{d I}~~ \frac{d J}{d V} \right] \left[
	\begin{array}{c}
	-\beta U(t) V(t)  \\
	\beta U(t)V(t)-\delta I(t)  \\
	pI(t) - cV(t)   
	\end{array}\right]\nonumber\\
	&=& \left[(1-\frac{U_s}{U(t)})~~1~~ \frac{\delta}{p} \right] \left[
	\begin{array}{c}
	-\beta U(t) V(t)  \\
	\beta U(t)V(t)-\delta I(t) \\
	pI(t) - cV(t)   
	\end{array}\right]\nonumber\\
	&=& (-\beta U(t) V(t) + U_s \beta V(t)) + (\beta U(t)V(t)-\delta I(t)) + (\delta I(t) - \frac{\delta c}{p} V(t))\nonumber\\
	&=& U_s\beta V(t) - \frac{\delta c}{p} V(t) = V(t) (U_s \beta  - \frac{\delta c}{p}).
	\end{eqnarray}
	Now, given $U_s \in [0,\Uc)$, with $\Uc = \frac{\delta c}{\beta p}$, it follows that $\dot{J}(x(t)) \leq 0$
	for every $x \in \X$ (note that it is not true that $\dot{J}(x(t)) < 0$ for $x \neq x_s$, as shown next, in Remark \ref{rem:nonas}). Then, $J$ is a Lyapunov function for system \eqref{eq:SysOrigAcut}, which means that each $x_s \in \setX_s^1$ is
	$\epsilon-\delta$ stable. Therefore, it is easy to see that the equilibrium set $\setX_s^1$ as a whole is also
	$\epsilon-\delta$ stable, which completes the proof. $\square$ 
\end{pf}

\begin{figure}
	\centering
	\includegraphics[width=0.65\textwidth]{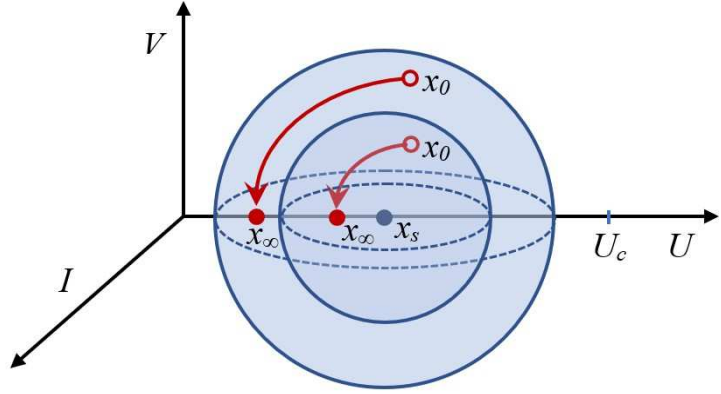}
	\caption{\small{Every point in $\setX_s^1$ is $\epsilon-\delta$ stable but not attractive. Initial states $x_0$ starting arbitrarily close
			to $x_s$ remains (for all $t \geq 0$) arbitrarily close to $x_s$, but does not converges to $x_s$. As a consequence, set $\setX_s^1$ is AS but the points
			inside it are not.}}
	\label{fig:SchemeStabil}
\end{figure}

\begin{rem} \label{rem:nonas}
	Note that, in the latter proof, it is not true that $\dot{J}(x(t)) < 0$ for every nonegative $x \neq x_s$. If for instance,
	the function $\dot{J}(x(t))$ is evaluated at $\hat x_s = (\hat U,0,0)$, with $\hat U \notin U_s$, we have that
    $\dot{J}(\hat x_s(t)) = 0$. In fact, $\dot{J}(x(t))$ is null along the whole $U$ axis, given that this axis is an
    equilibrium set. This means that the (individual) states in $\setX_s^1$ are $\epsilon-\delta$ stable, but not attractive.
\end{rem}
A schematic plot of such a behavior can be seen in Figure \ref{fig:SchemeStabil}.

\begin{rem} \label{rem:exaple2D}
	A similar behavior can be seen in system $\dot{x}=Ax$, when $A=[0~ -1;~0 -1~]$, or the 2-state Kermack-McKendrick epidemic model \cite{brauer2012mathematical,brauer2005kermack}: $\dot{S} = \beta SI$, $\dot{I}= \beta SI - \delta I$, being $S$ the susceptible and $I$ the infected individual. 
	In this latter model, $\RN:=(\delta/\beta)S_0$ and (the critical value for $S$ is) $S_c=\delta/\beta$. The AS set is given by all the states of the form $x_s:=(S_s,0)$, with $S_s \in [0,S_c)$. Furthermore, for this system, the maximum of $I$ occurs when $S=S_c$.
\end{rem}

%%%%%%%%%%%%%%%%%%%%%%%%%%%%%%%%%%%%%%%%%%%%%%%%%%%%%
\subsection{Asymptotic stability of $\setX_s^1$} \label{sec:ASproof}
%%%%%%%%%%%%%%%%%%%%%%%%%%%%%%%%%%%%%%%%%%%%%%%%%%%%%

In the next Theorem, based on the previous results concerning the attractivity and $\epsilon-\delta$ stability of $\setX_s^1$, the asymptotic
stability is formally stated.

\begin{theorem}\label{theo:AS}
	Consider system \eqref{eq:SysOrigAcut} constrained by the positive set $\mathbb X$.
	Then, the set $\setX_s^1$ defined in \eqref{ec:setXs1} is smallest asymptotically stable (AS) equilibrium set, with a domain of attraction given by $\setX$.
	%Furthermore, $\setX_s^1$ is the only AS set in $\setX$.
\end{theorem}

\begin{pf}
	The proof follows from Theorems \ref{theo:attract}, which states that $\setX_s^1$ is the smallest attractive in $\setX$, and \ref{theo:stab},
	which states the local $\epsilon-\delta$ stability of $\setX_s^1$. $\square$
\end{pf}

\vspace{0.5cm}

A critical consequence of the latter Theorem is that no equilibrium point in $\setX_s$ (neither in $\setX_s^1$, nor in $\setX_s^2$) can be used as setpoint in a control strategy design. The effect of antivirals (pharmocodynamic), for instance, is just to reduce the virus infectivity (by reducing 
the infection rate $\beta$) or the production of infectious virions (by reducing the replication rate $p$) \cite{hernandez2019modeling}. So, the previous stability analysis is still valid for controlled system, since only a modification of some of the parameters defining $\Uc$ is done. In such a context, only a controller able to consider the whole set $\setX_s^1$ as a target (a set-based control strategy, as zone MPC \cite{FerraLi10,gonz2020siMPC}) will be fully successful in controlling system \eqref{eq:SysOrigAcut}. 

%%%%%%%%%%%%%%%%%%%%%%%%%%%%%%%%%%%%%%%%%%%%%%%%%%%%%
\section{Characterization for different initial conditions} \label{sec:dyn_char}
%%%%%%%%%%%%%%%%%%%%%%%%%%%%%%%%%%%%%%%%%%%%%%%%%%%%%
In this section some further properties of system \eqref{eq:SysOrigAcut} concerning its dynamic are stated, based on the initial conditions at the infection time $t=0$. The objective is to fully characterize the states behavior in a qualitative way, including the times at which the virus and the infected cells reach their peaks. First, Property \ref{prop:Uinf} states some characteristics of $U_{\infty}$ for different initial conditions. Then, Theorem \ref{theo:key0}
states a general relationship between the peak times of $V$ and $I$ and the time at which $U$ reaches its critical value $\Uc$.

\begin{figure}
	\centering
	\includegraphics[width=0.65\textwidth]{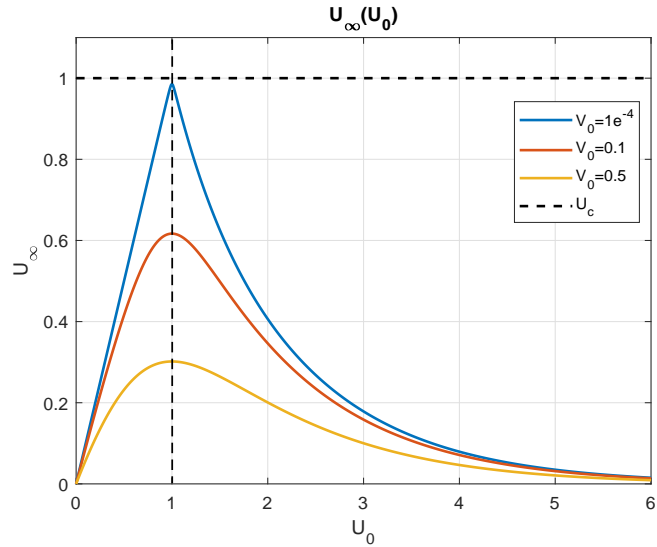}
	\caption{\small{According to equation \eqref{ec:dem10}, $U_\infty(U_0)$ is plotted for different values of $V_0$. 
	 All parameters are equal to $1$ for simplicity, which means that $\Uc=1$.}}
	\label{fig:UinftyU0}
\end{figure}

\begin{propt}\label{prop:Uinf}
	Consider system \eqref{eq:SysOrigAcut}, constrained by the positive set $\X$, at the beginning of the infection, i.e.,
	$U(0) = U_0 >0$, $I(0)= 0$ and $V(0)=V_0>0$ (i.e., $x(0)=(U(0),I(0),V(0)) \in \setX$).
	Consider also that $V_0$ is small enough. Then,
	\begin{enumerate}
		\item $U_\infty \to 0$ when $U_0\to\infty$ or $U_0\to 0$.
		\item $U_\infty \to \Uc$ when $U_0\to \Uc$.
		\item $0 < U_{\infty}( U_{0,1}, I_0, V_0) < U_{\infty}(U_{0,2},  I_0, V_0) < \Uc$, for initial conditions 
		$ U_{0,1} < U_{0,2} < \Uc$.
		\item $0 <  U_{\infty}( U_{0,2}, I_0, V_0) < U_{\infty}( U_{0,1}, I_0, V_0) < \Uc$, for initial conditions
		$\Uc <  U_{0,1} < U_{0,2}$. 
	\end{enumerate}
\end{propt}
\begin{pf}
	If $I_0=0$ and $V_0\approx 0$ then $\Ko\approx 0$. Therefore $W(-\RN e^{\Ko-\RN})\approx W(-\RN e^{-\RN})$, and  $U_\infty\approx -\Uc W(-\RN e^{-\RN})$  by \eqref{ec:dem10}.
	\begin{enumerate}
		\item $W(-\RN e^{-\RN })\to 0 $ when $-\RN e^{-\RN }\to 0$, which means that either $\RN\to 0$ or $\RN\to\infty$. This implies that $U_0\to 0$ 
		or $U_0\to\infty$, respectively.
		\item $W(-\RN e^{-\RN })\to -1$ when $-\RN e^{-\RN }\to -1/e$, which is true if $\RN \to 1$ or, the same, when $U_0\to \Uc$.
		\item Function $z(\RN )=\RN e^{-\RN }$ is strictly decreasing for $\RN \in (0,1)$ (note that ${\RN}_1 := \frac{c \delta U_{0,1}}{\beta p}$ and 
		${\RN}_2 := \frac{c \delta U_{0,2}}{\beta p}$ are in $(0,1)$, since they are smaller than $\Uc$), while $-W_p(\cdot)$ is strictly decreasing in $(-1/e,0)$. 
		So, $0 < -W_p(-{\RN}_1e^{-{\RN}_1}) < -W_p(-{\RN}_2e^{-{\RN}_2}) < 1$, which implies that 
		$0 < U_{\infty}( U_{0,1}, I_0, V_0) < U_{\infty}(U_{0,2},  I_0, V_0) < \Uc$.
		\item Function $z(\RN )=\RN e^{-\RN }$ is strictly increasing for $\RN \in (1,\infty)$, while $-W_p(\cdot)$ is strictly decreasing in $(-1/e,0)$.
		So, $0 < -W_p(-{\RN}_2e^{-{\RN}_2}) < -W_p(-{\RN}_1e^{-{\RN}_1}) < 1$, which implies that 
		$0 \textcolor{ao_english}{<} U_{\infty}( U_{0,2}, I_0, V_0) < U_{\infty}( U_{0,1}, I_0, V_0) < \Uc$. Figure \ref{fig:UinftyU0} shows $U_\infty$ as a function of $U_0$, taking $V_0$ as a parameter. $\square$
	\end{enumerate}
    
\end{pf}

\vspace{0.5cm}

\begin{figure}
	\centering
	\includegraphics[width=0.7\textwidth]{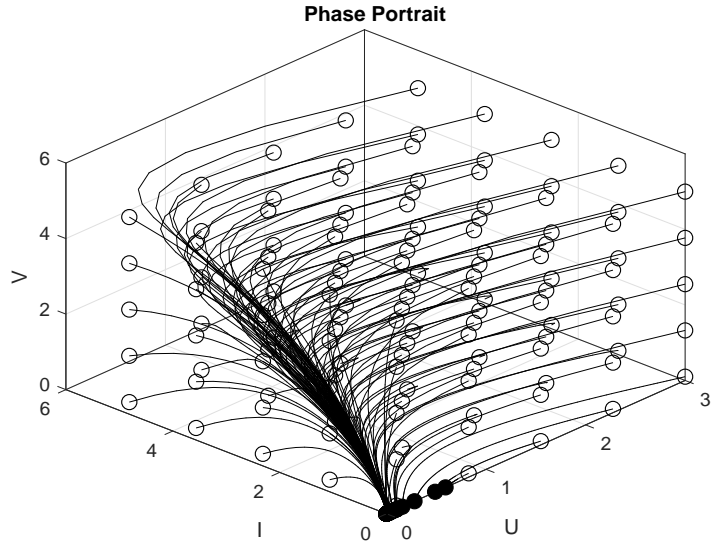}
	\caption{\small{Phase portrait of system \eqref{eq:SysOrigAcut}, with unitary parameters. Empty circles represent the initial state, while solid circles represent final states. Note that only the initial states with $U_0>\Uc=1$ corresponds to scenarios with $\RN>1$.}}
	\label{fig:PhaPor}
\end{figure}

Figure \ref{fig:PhaPor} shows a phase portrait of system \eqref{eq:SysOrigAcut}, where all parameters are equal to $1$ for simplicity, which means 
that $\Uc=1$.

\begin{theorem}[Virus behavior from the infection time]\label{theo:key0}
	Consider system \eqref{eq:SysOrigAcut}, constrained by the positive set $\X$, at the beginning of the infection, i.e.,
	$U(0) = U_0 >0$, $I(0)= 0$ and $V(0)=V_0>0$ (i.e., $x(0)=(U(0),I(0),V(0)) \in \setX$). If the virus spreads (according to Definition~\ref{defi:spread}), then there exist positive times $\check t_V$, $\hat t_I$, $t_c$ and $\hat t_V$, such that $\check t_V < \hat t_I < t_c < \hat t_V$, where $\check t_V$ and $\hat t_V$ are the times at which $V(t)$ reaches a local minimum and a local maximum, respectively, $\hat t_I$ is the time at which $I(t)$ reaches a local maximum, and $t_c$ is the time at which $U(t)$ reaches $\Uc$. Furthermore, $\dot V(t)<0$ for all $t> \hat t_V$.	
\end{theorem}

\begin{pf}
If $I(0)=0$ then $\Rn_V(0) < 1$, and Theorem~\ref{theo:key0} is a particular case of Theorem~\ref{theo:key} (ii).
\end{pf}

\vspace{0.5cm}

%Figure \ref{fig:UIV_max} shows a plot of the virus behavior.
%
%\begin{figure}
%	\centering
%	\includegraphics[width=0.75\textwidth]{Figuras/UIV_max.eps}
%	\caption{\small{Time evolution of $U$, $I$ and $V$, with unitary parameters $\beta,\delta,p,c$, for initial conditions $U_0=2,I_0=0,V_0=0.1$.}}
%	\label{fig:UIV_max}
%\end{figure}
%%

Next, a Remark concerning a particularity of Theorem \ref{theo:key0} is introduced, that may help to approximately determine the global maximum of the virus load.

\begin{rem}\label{cor:Vmax}
     Consider the hypothesis of Theorem \ref{theo:key0}. Then, it can be shown that for $V_0 \to 0$, $\hat t_V \to t_c$ from the right, and $\hat t_I \to t_c$ from the left, meaning that the peaks of $V$ and $I$ tends to occurs simultaneously at time $t_c$. This fact can be seen in Table \ref{tab:const}, when data coming from real patients are used to identify the model.
\end{rem}

A main consequence of Theorem \ref{theo:key0} and Remark \ref{cor:Vmax} is that it guaranties that the virus  will monotonically go to zero only after $U$ is below $U_c$. 
That is, any action devoted to steers $V$ to zero before $t_c$ may be counterproductive, since slowing down $V$ implies to soften the decreasing behavior of $U$, delaying time $t_c$, and maintaining $V$ large for a longer time.
This critical fact has a direct effect in a potential controlled system, when parameters $p$ or $\beta$ are scaled down by antiviral treatments. Assuming that $\RN$ remains greater than $1$, any attempt to steer $V$ to zero before $U$ has taken values below $\Uc$ may be unsuccessful, if 
$U$ is not controlled first.

%%%%%%%%%%%%%%%%%%%%%%%%%%%%%%%%%%%%%%%%%%%%%%%%%%%%%
\section{Particularization of the model with patient data} \label{sec:pac}
%%%%%%%%%%%%%%%%%%%%%%%%%%%%%%%%%%%%%%%%%%%%%%%%%%%%%
In this section, the parameters of model \eqref{eq:SysOrigAcut} will be associated to data from $9$ patients with COVID-19  - labeled as A,B,C,D,E,F,G,H and I - reported in \cite{woelfel2020clinical}. Here, we consider work in \cite{vargas2020host}, where different models for in-host SARS-CoV-2 were proposed and identified to fit the virus load data collected in \cite{woelfel2020clinical}. 

\begin{figure}
	\hspace{-1.5cm}
	\includegraphics[trim = 35mm 0mm 35mm 0mm, clip, width=1.2\textwidth]{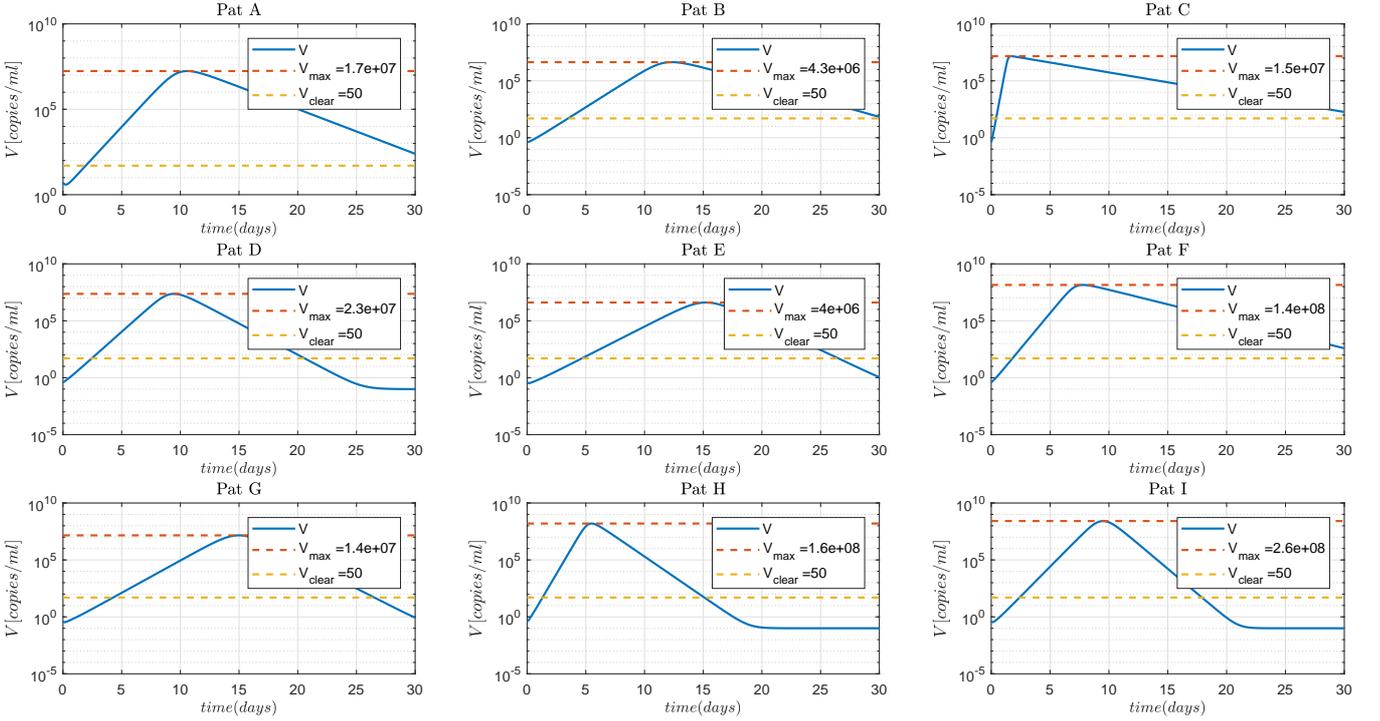}
	\caption{\small{Virus time evolution for all patients. As it can be seen, very different behaviors are obtained. $V_{clear}$ denotes
	a values of $50$ [copies/ml] under which the virus is considered cleared.}}
	\label{fig:V_all_pat}
\end{figure}

\begin{figure}
	\hspace{-1.5cm}
	\includegraphics[trim = 35mm 0mm 35mm 0mm, clip, width=1.2\textwidth]{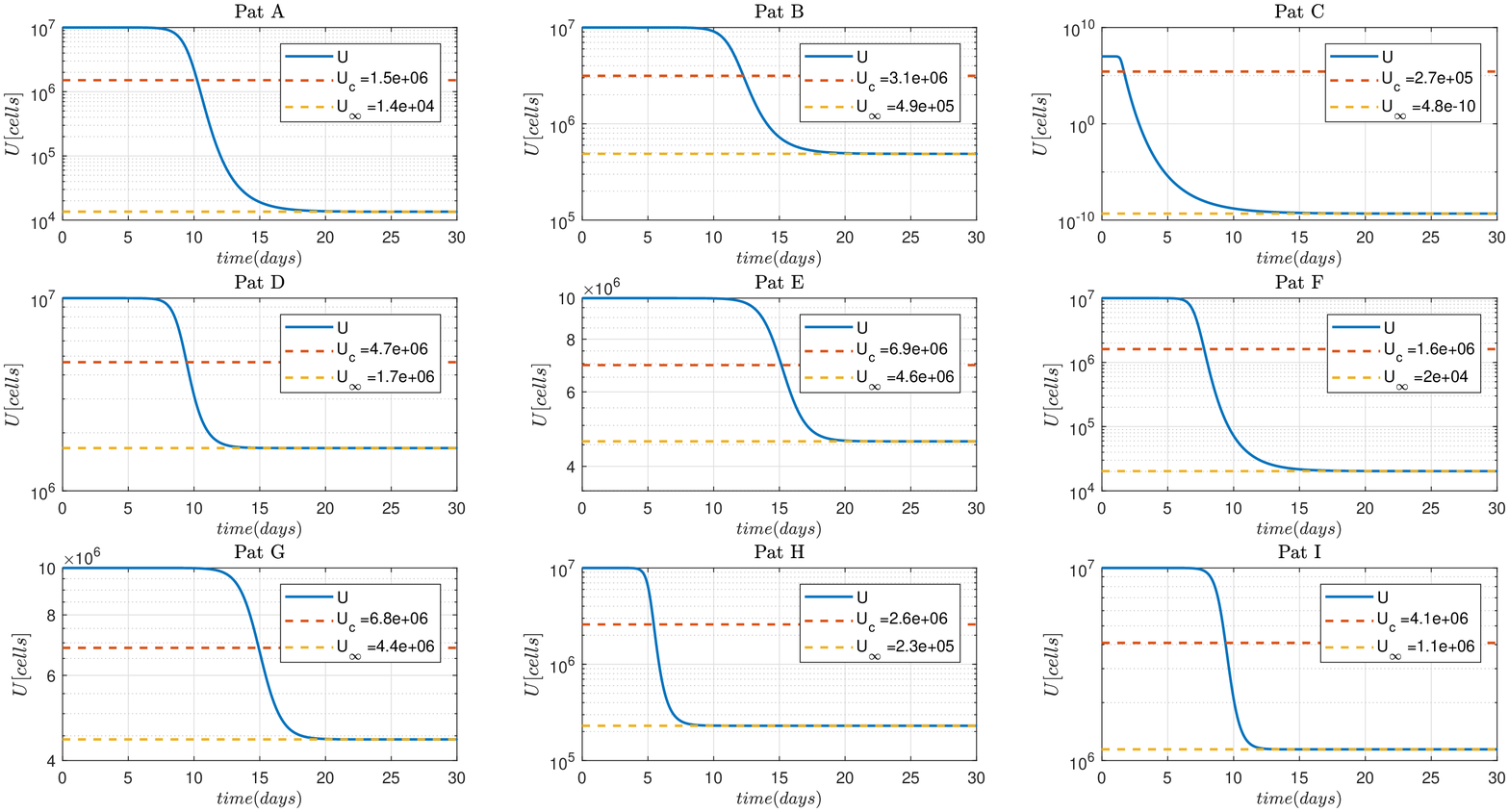}
	\caption{\small{Healthy cells time evolution for all patients. Pat C shows a very slow value of $U_{\infty}$ (practically zero), which
	suggest that the selected value of $U_0=1.0e^7$ is too large.}}
	\label{fig:U_all_pat}
\end{figure}

The initial number of target cells $U_0$ is estimated approximately $10^7$ cells \cite{vargas2020host}. $I_0$ is assumed to be $0$ while $V_0$ is determined by interpolation considering an incubation period of 7 days (note, that $V_0$ ranges from $0.02$ to $5.01$ $copies/mL$ which is below the detectable level of about $100$ $copies/mL$). Moreover, the onset of the symptoms is assumed to occurs 4 to 7 days after the infection time (day 0, Figure \ref{fig:V_all_pat} and \ref{fig:U_all_pat}).

Since the viral load is measured in logarithmic scale, the parameter fitting was performed minimizing the root mean square (RMS) difference on logarithmic scale between the model predictive output ($\hat{V}_i$) and the experimental measurements (${V}_i$), employing the Differential Evolution (DE) algorithm proposed in \cite{vargas2020host, hernandez2019modeling}. 
The parameters and the initial conditions ($U_0$, $I_0$ and $V_0$, with $t_0=0$ the infection time) of each patient are collected in Table \ref{tab:param}.

\begin{table}[H]
	\begin{center}		
		\caption{Target limited cell model parameter values for different patients with COVID-19 \cite{vargas2020host}.}
		\begin{tabular}{| c | c | c | c |c |}% c|c|c|}
			\hline
%			\multicolumn{5}{ |c| }{Amount of drug (mg)} \\ \hline
			Patient       &  $\beta$ &  $\delta$ &  $p$  & $c$ \\\hline%&  $U_0$ &  $I_0$ & $V_0$  \\ \hline
		A  &  $9.98\times10^{-8}$  &  0.61 &   9.3  &  2.3 \\%&	$1.0\times10^{7}$      &  0  &  5.01 \\%
		B  &  $1.77\times10^{-7}$  & 14.11 &  20.2  &  0.8\\%&	$1.0\times10^{7}$      &   0  &  0.31 \\
		C  &  $8.89\times10^{-7}$  & 79.51 & 134.4  &  0.4\\%&	$1.0\times10^{7}$      &   0  &  0.31	\\
		D  &  $3.15\times10^{-8}$  & 45.51 & 620.2  &  2.0\\%&	$1.0\times10^{7}$      &   0  &  0.31	\\
		E  &  $5.61\times10^{-8}$  &  7.51 &  96.4  &  5.0\\%&	$1.0\times10^{7}$      &   0  &  0.31	\\
		F  &  $1.41\times10^{-8}$  & 37.61 & 995.0  &  0.6\\%&	$1.0\times10^{7}$      &   0  &  0.31	\\
		G  &  $1.77\times10^{-8}$  &  8.21 & 338.4  &  5.0\\%&	$1.0\times10^{7}$      &   0  &  0.31	\\
	    H &	  $1.58\times10^{-8}$   & 21.11 & 927.8  &  1.8\\%&	$1.0\times10^{7}$      &   0  &  0.31	\\
        I &	  $4.46\times10^{-9}$  &  4.21 & 994.6   & 4.3\\\hline%&	$1.0\times10^{7}$      &   0  &  0.31	\\ \hline
		\end{tabular}
		\label{tab:param}
	\end{center}
\end{table}
According to the system analysis developed in previous sections, some relevant dynamical values are shown in Table \ref{tab:const}.
Constant $\alpha(0)$ (defined in Theorem \ref{theo:key}) is small than $10 \times 10^{-4}$ for all the patients, so it is not taken into account for the study.
\begin{table}[H]
	\begin{center}		
		\caption{Characterization Parameters of patients with COVID-19.}
		\begin{tabular}{| c | c | c | c | c | c | c | c | c | c | c | c |}
			\hline
			%			\multicolumn{5}{ |c| }{Amount of drug (mg)} \\ \hline
			Patient & $\Uc$ &$U_{\infty}$ & $\RN$ & $\Ko$ &  $\hat t_I$ &  $t_c$ & $\hat t_V$ & $V_{max}$ \\ \hline
			A	&	$1.51\times10^{6}$   & 	$1.36\times10^{4}$  & 6.61  & $-2.17\times10^{-7}$ & 10.16       &  10.24 &  10.58     &  $1.73\times10^{7}$   \\
			B	&	$3.15\times10^{6}$   & 	$4.88\times10^{5}$  & 3.18  & $-6.87\times10^{-8}$ & 11.54       &  12.26 &  12.32     &  $4.35\times10^{6}$   \\
			C	&	$2.66\times10^{5}$   & 	$4.81\times10^{-10}$ & 37.57 & $-6.89\times10^{-7}$          & 1.43        &  1.67  &  1.69      &  $1.47\times10^{7}$  \\
			D	&   $4.65\times10^{6}$   &  $1.67\times10^{6}$  & 2.15  & $-4.89\times10^{-9}$ & 9.04        &  9.42  &  9.44      &  $2.33\times10^{7}$ \\
			E	&   $6.94\times10^{6}$   &	$4.58\times10^{6}$  & 1.44  & $-3.48\times10^{-9}$ & 15.02       &  15.16 & 15.24      &  $4.03\times10^{6}$  \\
			F	&	$1.61\times10^{6}$   &	$2.03\times10^{4}$  & 6.21  & $-7.28\times10^{-9}$ & 7.12        &  7.76  &  7.78      &  $1.42\times10^{8}$  \\
			G	&	$6.84\times10^{6}$   &	$4.43\times10^{6}$  & 1.46  & $-1.1\times10^{-9}$ & 14.80       &  14.92 &  15.00       &  $1.44\times10^{7}$  \\
			H	&	$2.59\times10^{6}$   &	$2.3\times10^{5}$   & 3.86  & $-2.72\times10^{-9}$ & 5.16        &  5.44  &  5.48      &  $1.577\times10^{8}$  \\
			I	&	$4.08\times10^{6}$   & 	$1.14\times10^{6}$  & 2.45  & $-3.21\times10^{-10}$ & 9.28       &  9.38  &  9.50      &  $2.60\times10^{8}$   \\  \hline
		\end{tabular}
		\label{tab:const}
	\end{center}
\end{table}

Figures \ref{fig:V_all_pat} and \ref{fig:U_all_pat} shows the evolution of $V$ and $U$ for all patients. 
As expected, the states converges to $\setX_s^1$, although significantly different behavior can be observed for the different patients. From Figure \ref{fig:U_all_pat} it can be seen that the healthy cells final value $U_{\infty}$ is reduced in cases of patients with large values of $\RN$, in spite all patients have the same initial $U_0$. This can be explained from the fact that $W(\RN e^{-\RN}e^{\Ko})$ is monotonically decreasing for $\RN>1$ (see Figures \ref{fig:Lamb} and \ref{fig:FR1}), and therefore, $0<U_{\infty}(\mathcal{R}_{01})<U_{\infty}(\mathcal{R}_{02})$ for ${R}_{01}>{R}_{02}>1$. Note that the healthy cells of patient C converges to $U_{\infty}$ equals to $4.810\times10^{-10}~[cell]$, which can be explained by the fact that this patient has a reproduction number ($\RN$) of $37.57$, which is $5.2$ times above the cohort mean value of $7.21$. 
%It is important to note that the almost complete depletion of healthy cells of patient C must be understood in the context of susceptibility; that is, cells susceptible to virus infection and able to produce progeny virons. So, their depletion does not necessarily correspond to the decease of all cells of the respiratory tract. Indeed, it can be associated to the immune system inability to clear the infected cells and viral particles as well as the affinity of the viral particles to target cell receptors (i.e. angiotensin- converting enzyme 2, ACE2).  
Figure \ref{fig:V_all_pat} and Table \ref{tab:const} show that the viral load of patient C reaches the peak at $1.69$ days post infection (dpi) ($40.56$ hours post infection, hpi).% and, taking into account that the adaptive immune response takes in order of $7(dpi)$ to develop in COVID-19 patients \cite{woelfel2020clinical}, we can argue that in this case the viral evolution was mainly limited by susceptible cell availability. 

Furthermore, from Figure \ref{fig:V_all_pat}, it can be seen that for all the cases the viral load spreads (i.e.: the virus presents a peak)  although $\Rn_V(0)<0$ for all patients (i.e.,  $I_0=0$). This can be justified since $U_0\gg\Uc$ and, therefore, $\RN$ will be greater than $1+\alpha(0)$ for all patients (note that, $\alpha(0)<10 \times 10^{-4}$). Moreover, from Table \ref{tab:const}, we can corroborate that $\hat t_I>t_c>\hat t_V$ which is in accordance to what is stated in Theorem 4.1. %It is important to note that, in comparison with Influenza disease, which takes in order of $7-8 dpi$ to have a viral level below $V_{clear}$ level, the recovery period of the COVID-19 disease is in order of $25-30~dpi$. Therefore, considering that the alveolar epithelial cells (main target of SARS-CoV-2) proliferate in about 2-3 weeks, the acute assumption of the model may be reconsidered, by adding some kind of target cell production and natural death rate (apoptosis) terms in equation (\ref{eq:SysOrigAcut}.a) \cite{ciupe2017host}. 

Concerning the immune response, this model makes the assumption that it is constant and independent on viral load as well as infected cells. Furthermore, neither innate or adaptive response are modeled, being the viral load dynamic mainly limited by target cells availability. Since recent studies have shown a dysfunctional immune response (i.e.: lymphogenia, desregulated secretion of pro-inflammatory cytokines, excessive infiltration of monocytes, macrophages and T cells, among others) \cite{tay2020trinity,diao2020reduction}, this effect should be added in the proposed model, in order to have a more reliable representation (and, eventually, a more realistic control objective). In addition, a more reliable standard to measure the severity of disease could be related with the viral spreadability as well as the deregulated inflammatory response.

%%
%\begin{figure}
%	\centering
%	\includegraphics[width=0.65\textwidth]{Figuras/Virus_all_pat.eps}
%	\caption{\small{Virus time evolution for all patients. As it can be seen, very different behaviors are obtained}}
%	\label{fig:V_all_pat}
%\end{figure}
%%
%
%
%%
%\begin{figure}
%	\centering
%	\includegraphics[width=0.65\textwidth]{Figuras/HeCells_all_pat.eps}
%	\caption{\small{Healthy cells time evolution for all patients. Pat C shows a very slow value of $U_{\infty}$ (practically zero), which
%	suggest that the selected value of $U_0=1.0e^7$ is too large.}}
%	\label{fig:U_all_pat}
%\end{figure}
%%
%
%

%%%%%%%%%%%%%%%%%%%%%%%%%%%%%%%%%%%%%%%%%%%%%%%%%%%%%
\section{Conclusions}\label{sec:conc}
%%%%%%%%%%%%%%%%%%%%%%%%%%%%%%%%%%%%%%%%%%%%%%%%%%%%%
In this work a full dynamical characterization of a COVID-19 in-the-host target-cell model is performed. Opposite to what happens in other similar models, it is shown that there exists a minimal non-punctual stable equilbrium set depending only on the system parameters. Furthermore, it is shown that there exists a parameter-depending threshold for the susceptible cells that fully characterizes the virus and infected cells qualitative behavior. Simulations performed with real-patient data demonstrate the potential utility of such system dynamic characterization to tailor the the most valuable pipeline drugs against SARS-CoV-2.

%%%%%%%%%%%%%%%%%%%%%%%%%%%%%%%%%%%%%%%%%%%%%%%%%%%%%%%%%%%%%%%%%%%%%%%%%%%%%%%%%%%%%%%%%%%%%%%%
\section{Appendix 1. Stability theory}\label{sec:app1}
%%%%%%%%%%%%%%%%%%%%%%%%%%%%%%%%%%%%%%%%%%%%%%%%%%%%%%%%%%%%%%%%%%%%%%%%%%%%%%%%%%%%%%%%%%%%%%%%

In this section some basic definitions and results are given concerning the asymptotic stability of sets and Lyapunov theory,
in the context of non linear continuous-time systems. All the following definitions are referred to system 
\begin{eqnarray}\label{eq:difeqini}
\dot x (t) = f(x(t)),~~x(0)=x_0,
\end{eqnarray}
where $x$ is the system state constrained to be in $\X \subseteq \R^n$, $f$ is a Lipschitz continuous nonlinear function,
and $\phi(t;x)$ is the solution for time $t$ and initial condition $x$.
\begin{defi}[Equilibrium set]
	Consider system \ref{eq:difeqini} constrained by $\X$. The set $\setX_s \subset \X$ is an equilibrium set if each point $x \in \setX_s$
	is such that $f(x)=0$ (this implying that $\phi(t;x)=x$ for all $t \geq 0$).
\end{defi}
\begin{defi}[Attractivity of an equilibrium set]\label{def:attrac_set}
	Consider system \ref{eq:difeqini} constrained by $\X$. A closed equilibrium set $\setX_s \subset \X$ is attractive in $\setX \subset \X$
	if $\lim_{t \rightarrow \infty} \|\phi(t;x)\|_{\setX_s} =0$ for all $x \in \setX$.
\end{defi}
Any set containing an attractive set is attractive, so the significant attractive of a constrained system set is given by the smallest one.
\begin{defi}[$\epsilon-\delta$ local stability of an equilibrium set]\label{def:eps_del_stab}
	Consider system \ref{eq:difeqini} constrained by $\X$. A closed equilibrium set $\setX_s \subset \X$ is $\epsilon-\delta$ locally stable
	if for all $\epsilon >0$ it there exists $\delta>0$ such that in a given boundary of $\setX_s$, $\|x\|_{\setX_s} <\delta$, it
	follows that $ \|\phi(t;x)\|_{\setX_s} < \epsilon$, for all $t \geq 0$.
\end{defi}
\begin{defi}[Asymptotic stability (AS) of an equilibrium set]\label{def:AS}
	Consider system \ref{eq:difeqini} constrained by $\X$. A closed equilibrium set $\setX_s \in \X$ is asymptotically stable (AS) in $\setX \subset \X$
	if it is $\epsilon-\delta$ locally stable and attractive in $\setX$.
\end{defi}
\begin{theorem}[Lyapunov theorem]
	Consider system \ref{eq:difeqini} constrained by $\X$ and an equilibrium state $x_s \in \setX_s \subset \X$. Let consider a 
	function $V(x): \R^n \rightarrow \R$ such that $V(x)>0$ for $x \neq x_s$, $V(x_s)=0$ and $\dot{V}(x(t)) \leq 0$, denoted as Lyapunov function.
	Then, the existence of such a function implies that $x_s \in \setX_s$ is $\epsilon-\delta$ locally stable. If in addition
	$\dot{V}(x(t)) < 0$ for all $x \neq x_s$ and $\dot{V}(x_s) = 0$, then $x_s \in \setX_s$ is asymptotically stable.
\end{theorem}
%
%The next definition is referred to the controlled system 
%%
%\begin{eqnarray}\label{eq:difeqinicont}
%\dot x (t) = f(x(t),\eta(t)),~~x(0)=x_0,
%\end{eqnarray}
%%
%where $\eta(t) \in \U \subset \R^m$ represents de control actions, and $\phi(t;x,\eta(\cdot))$ is the solution of the system for time $t$, 
%initial condition $x$ and a control function $\eta(\cdot):[0,t] \to \R^m$.
%%
%\begin{defi}[Controlled equilibrium]\label{def:cont_set}
%	Let consider system \ref{eq:difeqinicont} constrained by $\X$. A closed set $\setX_s \in \X$ is a controlled 
%	equilibrium set if it there exists an fixed input $\eta \in \U$ such that $f(x,\eta) = 0$ for every $x \in \setX_s$.
%\end{defi}
%%
%%
%\begin{defi}[Stabilizable set]\label{def:stab_set}
%	Let consider system \ref{eq:difeqinicont} constrained by $\X$. A closed controlled equilibrium set $\setX_s \in \X$ is stabilizable 
%	w.r.t. $\setX$ if every state in $\setX$ can be steered and maintained in $\setX_s$ by means of feasible control actions $\eta(\cdot) \in \U$
%    and following feasible trajectories, $\phi(t;x,\eta(\cdot)) \in \X$ for all $t\geq0$. Set $\setX$ is then denoted as stabilizable set.
%\end{defi}
%%
%Clearly, the stabilizable set definition accounts for the AS of an equilibrium set, but considering the existence of control actions and
%the feasibility of both, states and inputs along the path.

%%%%%%%%%%%%%%%%%%%%%%%%%%%%%%%%%%%%%%%%%%%%%%%%%%%%%%%%%%%%%%%%%%%%%%%%%%%%%%%%%%%%%%%%%%%%%%%%
\section{Appendix 2. Derivation of the basic reproduction number $\RN$}\label{sec:app2}
%%%%%%%%%%%%%%%%%%%%%%%%%%%%%%%%%%%%%%%%%%%%%%%%%%%%%%%%%%%%%%%%%%%%%%%%%%%%%%%%%%%%%%%%%%%%%%%%
The derivation of the basic reproduction number $\RN$ will be given by means of the concept of next-generation 
matrix \cite{van2017reproduction}. Consider system \eqref{eq:SysOrigAcut} and assume that a healthy equilibrium exists, of the form $x_0=(U_0,0,0)$, and it is stable in absence of disease.
Of the complete state of system \eqref{eq:SysOrigAcut}, $x=(U,I,V)$, only two states depend 
on infected cells, that is $I$ and $V$. Let us rewrite the ODEs for this two states in the form
\begin{eqnarray*}
\dot{I}(t)&=&\mathcal{F}_I(x)-\mathcal{G}_I(x)\\
\dot{V}(t)&=&\mathcal{F}_V(x)-\mathcal{G}_V(x)
\end{eqnarray*}
where $\mathcal{F}_i(x)$, $i=\{I,V\}$, is the rate of appearance of new infections in compartment $i$, while $\mathcal{G}_i(x)$, $i=\{I,V\}$, is the rate of other transitions between compartment $i$ and the other
infected compartments, that is
\[
\begin{array}{lll}
\mathcal{F}_I(x)=\beta U(t)V(t) & \mbox{ and } & \mathcal{G}_I(x)= \delta I(t)\\
\mathcal{F}_V(x)= 0 & \mbox{ and }& \mathcal{G}_V(x)= -p I(t) +c V(t)
\end{array}
\]

If we now define
\[
F=  \bmat{cc} \displaystyle \frac{\partial \mathcal{F}_I(x)}{\partial I} & \displaystyle \frac{\partial \mathcal{F}_I(x)}{\partial V}\\\\
\displaystyle \frac{\partial \mathcal{F}_V(x)}{\partial I} & \displaystyle \frac{\partial \mathcal{F}_V(x)}{\partial V}\emat_{x=x_0} = \bmat{cc} 0 & \beta U_0 \\ 0 & 0\emat
\]
and
\[
G= \bmat{cc} \displaystyle \frac{\partial \mathcal{G}_I(x)}{\partial I} & \displaystyle \frac{\partial \mathcal{G}_I(x)}{\partial V}\\\\
\frac{\displaystyle \partial \mathcal{G}_V(x)}{\partial I} & \displaystyle \frac{\partial \mathcal{G}_V(x)}{\partial V}\emat_{x=x_0} = \bmat{cc} \delta & 0 \\ -p & c\emat
\]
then matrix $FG^{-1}$, represents the so-called \textit{next-generation matrix}. Each $(i,j)$ entry of such a matrix represents the expected number of secondary infections in compartment $i$ produced by an infected cell introduced in compartment $j$. The spectral radius of this matrix, that is, the maximum absolute value of its eigenvalues, defines the basic reproduction number $\RN$. 

For the specific case of system \eqref{eq:SysOrigAcut}, the \textit{next-generation matrix} is given by
\[
FG^{-1}= \bmat{cc} \displaystyle \frac{\beta p U_0}{c \delta} & \displaystyle \frac{\beta U_0}{c}\\ \\ 0 & 0\emat
\]
Therefore, the basic reproduction number $\RN$ is given by
\[
\RN=:\frac{\beta p U_0}{c \delta} 
\]

Notice that $\RN$ coincides with the entry $(1,1)$ of matrix $FG^{-1}$, thus meaning that $\RN$ represents the expected number of secondary infections produced in compartment $I$ by an infected cell originally in $I$.

%%%%%%%%%%%%%%%%%%%%%%%%%%%%%%%%%%%%%%%%%%%%%%%%%%%%%%%%%%%%%%%%%%%%%%%%%%%%%%%%%%%%%%%%%%%%%%%%
\section{Appendix 3. General virus characterization for the dynamic \eqref{eq:SysOrigAcut}} \label{sec:app3}
%%%%%%%%%%%%%%%%%%%%%%%%%%%%%%%%%%%%%%%%%%%%%%%%%%%%%%%%%%%%%%%%%%%%%%%%%%%%%%%%%%%%%%%%%%%%%%%%

The next theorem characterizes all possible virus behavior, depending on the arbitrary initial conditions and the parameters. 
\begin{theorem}[Virus behavior from an arbitrary time]\label{theo:key} 
	Consider system \eqref{eq:SysOrigAcut}, constrained by the positive set $\X$, at some arbitrary time $t_0$, with
	$U(t_0) >0$, $I(t_0) \geq 0$ and $V(t_0) > 0$ (i.e., $x(t_0)=(U(t_0),I(t_0),V(t_0)) \in \setX$).
	Let define $\Rn(t):=\frac{\beta p U(t)}{\delta c}$ and $\Rn_V(t):=\frac{p I(t)}{c V(t)}$. Then,
	\begin{enumerate}
		\item if $\Rn_V(t_0) < 1$ and $\Rn(t_0) < 1 + \alpha(t_0)$, where $\alpha(t_0)$ is a positive value, depending on $I(t_0),V(t_0)$ and the parameters,
		 then, $\dot V(t)<0$ for all $t > t_0$ (i.e., $V(t)$ is strictly decreasing for all $t > t_0$ and, so, the virus does not spread in the body host),
		\item if $\Rn_V(t_0) < 1$ and the virus spreads in the body host (i.e. $V(t)$ reaches a local maximum, at some time $\hat t_V > t_0$,
		before it goes to zero), then $\Rn(t_0) > 1 + \alpha(t_0)$. Even more, $V(t)$ has a local minimum at some time $\check t_V>t_0$,
		with $\check t_V<\hat t_V$, and after the 
		local maximum it is strictly decreasing, i.e., $\dot V (t)<0$ for all $t > \hat t_V$. On the other hand, $I(t)$ has only one local maximum
		at time $\hat t_I$, being $t_0 <\check t_V < \hat t_I < t_c < \hat t_V$, where $t_c$ is the time when $\Rn$ reaches $1$ from above
		(or, the same, when $U(t)$ reaches its critical value $\setU_c$).
		\item if $\Rn_V(t_0) > 1$, then it there exists $\hat t_V > t_0$ such 
		that $\dot V(t)>0$ for all $t_0 < t < \hat t_V$, while $\dot V(t)<0$ for all $t > \hat t_V$ (i.e., $V(t)$ has 
		a global maximum at $\hat t_V$ and then it is strictly decreasing for all $t > \hat t_V$). This means that 
		the virus spreads in the host body.
	\end{enumerate}
\end{theorem}
\begin{pf}
	First note that for positive parameters, equation (\ref{eq:SysOrigAcut}.a) implies that $\dot U(t)<0$ for all $t \geq t_0$ and, so, 
	$\Rn(t)=U(t)\frac{\beta p}{c \delta}$ is strictly decreasing for all $t\geq t_0$. Given that 
	$\Rn_{\infty}:=\lim_{t \to \infty}\Rn(t)=U_{\infty}\frac{\beta p}{c \delta}$ and $U_{\infty} < \Uc := \frac{c \delta}{\beta p}$ (see 
	Theorem \ref{theo:attract}), then $\Rn_{\infty}<1$. So, for $\Rn(t_0)>1$, there exists only one time 
	\begin{eqnarray}\label{eq:theo0}
	t_c>t_0,
	\end{eqnarray}
	at which $\Rn(t_c)=1$, being $\Rn(t)>1$ for $t_0<t<t_c$ and $\Rn(t)<1$ for $t>t_c$.
	\begin{enumerate}
		\item By hypothesis, $\Rn_V(t_0) = \frac{p I(t_0)}{c V(t_0)} < 1$, which implies that $\dot V(t_0) < 0$ ($V$ 
		starts decreasing at $t_0$). On the other hand, by Lemma \ref{lem:extremes}, $V$ reaches a minimum or an inflection point at 
		at time $t_V^*>t_0$ if $\Rn(t_V^*)\geq 1$. But $\Rn(t)$ is strictly decreasing and, so, to reach 1 (at least) at $t_V^* > t_0$ it must be
		\begin{eqnarray}\label{eq:theo1}
		\Rn(t_0) = 1 + \alpha(t_0),
		\end{eqnarray}
		for some $\alpha(t_0)>0$. However, by hypothesis it is $\Rn(t_0) < 1 + \alpha(t_0)$, which implies that $\dot V(t)<0$ for all $t > t_0$.
		This concludes the proof.
		\item By hypothesis, $\Rn_V(t_0)<1$, which means that $\dot V(t_0) < 0$, and $V(t)$ reaches a local maximum at some time $\hat t_V > t_0$.
		Therefore, $V(t)$ must reach a local minimum at some time $t_0 < \check t_V < \hat t_V$.
		At both, the local minimum and maximum it is $\dot V(\check t_V)=0$ and $\dot V(\hat t_V)=0$. Then, by Lemma \ref{lem:extremes} it is 
		$\Rn(\check t_V)>1$ and $\Rn(\hat t_V)<1$, respectively. 
		
		Since $\Rn(t)$ is strictly decreasing for $t>t_0$ and $\check t_V >t_0$ (note that $\dot V(t_0)<0$ and $\dot V(\check t_V)=0$, so $\check t_V$
		cannot be equal to $t_0$.), then $\Rn(t_0)>\Rn(\check t_V)>1$, which implies that $\Rn(t_0) > 1 + \alpha(t_0)$ for some $\alpha(t_0)>0$.
		Furthermore, $t_0<\check t_V < t_c < \hat t_V$ ($\Rn$ crosses $1$ at time $t_c$ between $\check t_V$ and $ \hat t_V$) and, given that $\Rn(\hat t_V)<1$, $V$ cannot reaches another local minimum after its local maximum.
		This implies that $\dot V(t)<0$ for $t > \hat t_V$.
		
		From the minimum and maximum conditions of $V$, at times $\check t_V$ and $\hat t_V$, they are $\dot V(\check t_V)=0$, $\ddot{V}(\check t_V)>0$ 
		and $\dot V(\hat t_V)=0$, $\ddot{V}(\hat t_V)<0$, respectively. After some algebraic computation, it is easy to see that 
		$\dot I(\check t_V) > 0$ and $\dot I(\hat t_V) <0$, which means that $I(t)$ must reach a maximum at some time $\hat t_I$, fulfilling
		$\check t_V < \hat t_I < \hat t_V$. Even more, it must be $\dot I(\hat t_I) =0$, or
		\begin{eqnarray}\label{eq:theo2}
		\beta U(\hat t_I) V(\hat t_I)-\delta I(\hat t_I)=0.
		\end{eqnarray}
		Given that $\dot V(t) >0$ for $\check t_V < t < \hat t_V$ (it goes from its minimum to its maximum), then 
		by (\ref{eq:SysOrigAcut}.a), $I(\hat t_I) > \frac{c}{p}V(\hat t_I)$. Replacing this later condition in \eqref{eq:theo2},
		it follows that
		\begin{eqnarray}\label{eq:theo3}
		(\beta U(\hat t_I) - \frac{\delta c}{p}) V(\hat t_I)> \beta U(\hat t_I) V(\hat t_I)-\delta I(\hat t_I) = 0,
		\end{eqnarray}
		which implies that $\Rn(\hat t_I) = \frac{\beta p U(\hat t_I)}{\delta c} >1$ and, then, $\hat t_I < t_c$. 
		Therefore, $t_0 < \check t_V  < \hat t_I < t_c < \hat t_V$, which concludes the proof.
		\item By hypothesis, $\Rn_V(t_0) = \frac{p I(t_0)}{c V(t_0)} > 1$, which implies $\dot V(t_0) = pI(t_0)-cV(t_0) > 0$ 
		($V$ starts increasing at $t_0$). Since $V_{\infty}=0$ (Theorem \ref{theo:attract}),	then, there exists $\hat t_V> t_0$ such that 
		$V(\hat t_V)$ is a maximum. According to Lemma \ref{lem:extremes}, if $V$ has a maximum at $\hat t_V$, then $\Rn(\hat t_V) <1$.
		
		On the other hand, for $V(t)$ to reach a minimum after time $\hat t_V$, it must be $\Rn(\hat t_V)>1$. But $\Rn(t)$ is strictly decreasing
		for $t>t_0$, which means that no further minimum exists after $\hat t_V$. This implies that $\dot V(t)<0$ for all $t >\hat t_V$,
		which concludes the proof.
		$\square$
	\end{enumerate}	
\end{pf}

\vspace{0.5cm}

\begin{lem}\label{lem:extremes}
	Consider system \eqref{eq:SysOrigAcut}, constrained by the positive set $\X$, at some arbitrary time $t_0$, with
	$U(t_0) >0$, $I(t_0) \geq 0$ and $V(t_0) > 0$ (i.e., $x(t_0)=(U(t_0),I(t_0),V(t_0)) \in \setX$). 
	%Consider that at some time $t_V^*$, it is $\dot V(t_V^*)=0$. 
	Then, 
	(i) if $V(t)$ reaches a local minimum at time $t_V^* >t_0$, then $\Rn(t_V^*) > 1$,
	(ii) if $V(t)$ reaches a local maximum at time $t_V^*>t_0 $, then $\Rn(t_V^*) < 1$, and
	(iii) if $V(t)$ reaches an inflection point at time $t_V^*>t_0$ (a point in which $\dot V=0$ and $\ddot{V}=0$), then
	$t_V^*=t_c$, where $t_c$ is the (unique) time at which $\Rn$ reaches 1 (i.e., $\Rn(t_c) = 1$ or, the same, $U(t_c)=\Uc$).
\end{lem}

\begin{pf}
	Any of the three hypothesis ($V(t)$ reaches a local minimum, a local maximum or a inflection point) implies that  
	\begin{eqnarray}\label{eq:theoK0}
		\dot{V}(t_V^*) = pI(t_V^*) - c V(t_V^*) = 0,
	\end{eqnarray}
	which means that 
	\begin{eqnarray}\label{eq:theoK0_5}
	I(t_V^*) = p/c I(t_V^*).
	\end{eqnarray}
	Consider the critical case of an inflection point, i.e., 
	\begin{eqnarray}\label{eq:theoK1}
	\ddot{V}(t_V^*) = p\dot I(t_V^*) - c \dot V(t_V^*) = p \dot I(t_V^*) = 0.
	\end{eqnarray}
	From \eqref{eq:theoK1} it is $\dot I(t_V^*) = 0$ which, by (\ref{eq:SysOrigAcut}.b) at $t_V^*$, is equivalent to 
	\begin{eqnarray}\label{eq:theoK2}
	\dot I(t_V^*) = \beta U(t_V^*) V(t_V^*) - \delta I(t_V^*) = 0.
	\end{eqnarray}
	Now, by \eqref{eq:theoK0_5}, we have
	\begin{eqnarray}\label{eq:theoK3}
	(\frac{\beta p}{c} U(t_V^*) - \delta ) I(t_V^*) = 0.
	\end{eqnarray}
	Given that $I(t_V^*) > 0$ (note that $I(t)$ is positive for all $t>0$), then $\frac{\beta p}{c} U(t_V^*) - \delta = 0$, or
	\begin{eqnarray}\label{eq:theoK4}
	\Rn(t_V^*) = \frac{\beta p}{c \delta} U(t_V^*) = 1.
	\end{eqnarray}
	This way if an inflection point does occurs at $t_V^*$, then $t_V^*=t_c$, where $t_c$ is the time at which 
	$\Rn=1$. This proves item (iii).

	Furthermore, if $V$ reaches a local minimum at $t_V^*$, then $\ddot{V}(t_V^*)>0$ (instead of $\ddot{V}(t_V^*)=0$,
	as it is in (\ref{eq:theoK1}), which by \eqref{eq:theoK0_5} implies that
	\begin{eqnarray}\label{eq:theoK5}
	\Rn(t_V^*) = \frac{\beta p}{c \delta} U(t_V^*) > 1.
	\end{eqnarray}
	This proves item (i).
	
	On the other hand, if $V$ reaches a local maximum at $t_V^*$, then $\ddot{V}(t_V^*)<0$ (instead of $\ddot{V}(t_V^*)=0$,
	as it is in (\ref{eq:theoK1}), which by \eqref{eq:theoK0_5} implies that
	\begin{eqnarray}\label{eq:theoK6}
	\Rn(t_V^*) = \frac{\beta p}{c \delta} U(t_V^*) < 1.
	\end{eqnarray}
	This proves item (ii). $\square$
\end{pf}

\vspace{0.5cm}

\begin{figure}
	\centering
	\includegraphics[width=0.9\textwidth]{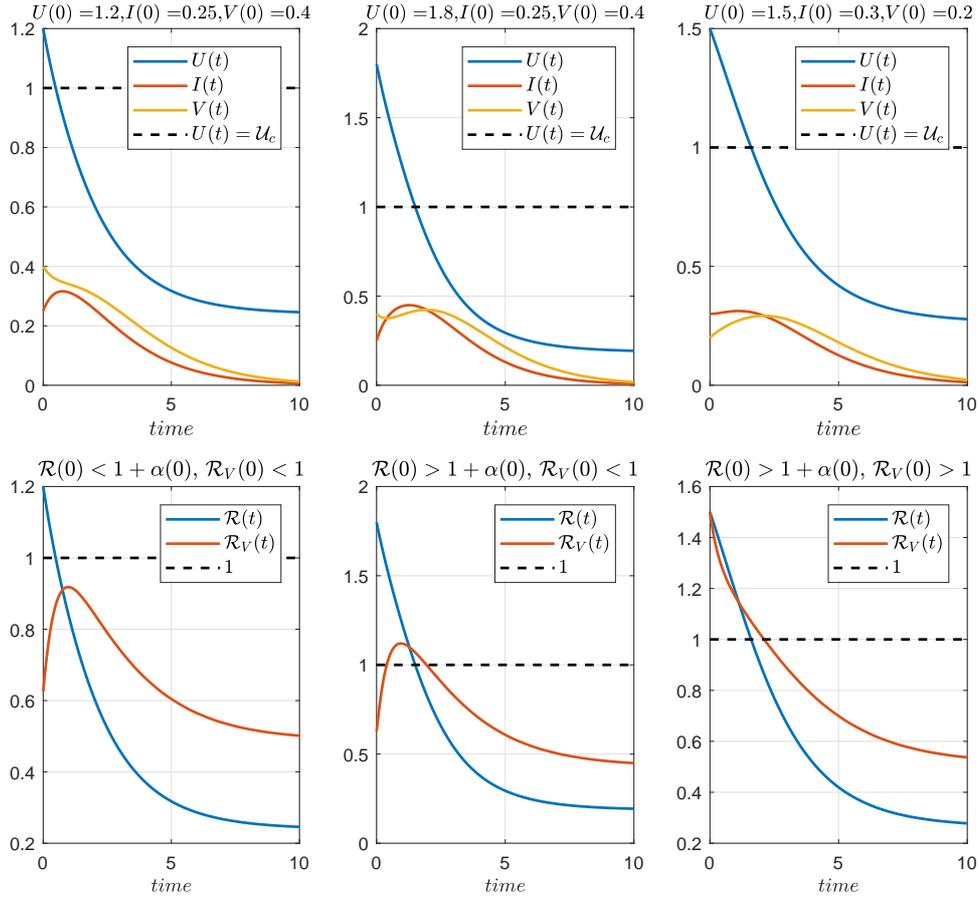}
	\caption{\small{Theorem \ref{theo:key}, schematic plot. $\beta=\delta=p=c=1$, which means that $\Uc=1$. The initial time is selected to be 
			$t_0=0$. First column: case (i). System time evolution when $U(0)=1.2$, $I(0)=0.25$ and $V(0)=0.4$, such that $\Rn_V(0)=0.63 <1$ and $\Rn(0)=1.2 < 1+\alpha(0)$, being $\alpha(0)=0.43$. Second column: case (ii). System time evolution when $U(0)=1.8$, $I(0)=0.25$ and $V(0)=0.4$, such 
			that $\Rn_V(0)=0.63<1$ and $\Rn(0)=1.8> 1+\alpha(0)$, being $\alpha(0)=0.43$. Third column: case(iii). System time evolution when $\Rn_V(0)>1$ and $\Rn(0)> 1+\alpha(0)$.}}
	\label{fig:main}
\end{figure}

Figure \ref{fig:main} shows schematic plots of cases i, ii and iii of Theorem \ref{theo:key}. For the sake of simplicity, $\beta=\delta=p=c=1$, 
which means that $\Uc=1$, while the initial time is selected to be zero, i.e., $t_0=0$.

Figure \ref{fig:main}, first column, illustrates case (i), in which $\Rn_V(0)<1$ and $\Rn(0) < 1+\alpha(0)$, being $\alpha(0)=0.43$ (numerically computed). 
This simulation shows that even if $\Rn(0) > 1$, but it is not greater than $1+\alpha(0)$, the virus does not spread in the body host, i.e, it is strictly decreasing. Note that the maximum of $\Rn_V$ does not reach $1$.
Figure \ref{fig:main}, second column, illustrates case (ii), in which $\Rn_V(0)<1$ and $\Rn(0)> 1+\alpha(0)$, being $\alpha(0)=0.43$. 
This simulation shows that for $\Rn(0) > 1+\alpha(0)$, the maximum of $\Rn_V$ is greater than $1$, and the period of time in which $\Rn_V>1$ 
is precisely the period of time the virus increases and reach a maximum (i.e., it spreads in the body host). As expected, $V$ reaches a
minimum first, then $I$ reaches a maximum, then $\Rn$ reaches $1$ (at $t_c$) and, finally, $V$ reaches a maximum, before to strictly decrease
to zero.
Figure \ref{fig:main}, third column, shows case(iii), in which $\Rn_V(0)>1$ and $\Rn(0)> 1+\alpha(0)$. As expected, the virus 
spreads on the body host an it has only one global maximum. 
%
%%
%\begin{figure}
%	\centering
%	\includegraphics[width=0.65\textwidth]{Figuras/main2.eps}
%	\caption{\small{Theorem \ref{theo:key}, }}
%	\label{fig:main2}
%\end{figure}
%%

\begin{rem}\label{rem:alphaneg}
	The value of $\alpha$ is necessary to properly understand and characterize the system behavior according to the initial conditions and parameters.
	Although it cannot be explicitly defined, it can be computed numerically. Furthermore, it should be noted that for real patients data, this values 
	use to be small in comparison with $\Rn$, given that $\beta$ is	small (see Table \ref{tab:const}). 
%	Therefore, its effect will be neglected for the real patient data simulations (Table \ref{tab:alfas} shows the value of $\alpha(t_{tr})$
%	when the treatment starts at day seven from the infection).
\end{rem}

%----------------------------------------------------------------------------------------
%	BIBLIOGRAPHY
%----------------------------------------------------------------------------------------

\bibliographystyle{elsarticle-num} 
\bibliography{biblioIR}

\end{document}